\documentclass[reqno,10pt,centertags]{amsart}
\usepackage{amssymb}
\usepackage{mathptmx}
\usepackage{hyperref}
\DeclareSymbolFont{SY}{U}{psy}{m}{n}
\DeclareMathSymbol{\emptyset}{\mathord}{SY}{'306}

\renewcommand{\eqref}[1]{{\rm(\ref{#1})}}
\newcommand{\Ref}[1]{{\rm\ref{#1}}}

\newcommand{\bbC}{{\mathbb C}}
\newcommand{\bbR}{{\mathbb R}}

\newcommand{\cA}{{\mathcal A}}
\newcommand{\cB}{{\mathcal B}}

\newcommand{\cH}{{\mathcal H}}

\newcommand{\cK}{{\mathcal K}}

\newcommand{\cO}{{\mathcal O}}
\newcommand{\cP}{{\mathcal P}}

\newcommand{\cW}{{\mathcal W}}

\newcommand{\hsE}{\widehat{\sf E}}

\newcommand{\ri}{{\rm i}}

\newcommand{\sE}{{\sf E}}

\newcommand{\no}{\nonumber}

\newcommand{\tA}{\widetilde{A}}

\newcommand{\fH}{\mathfrak{H}}
\newcommand{\fK}{\mathfrak{K}}

\newcommand{\Max}{\mathop{{\rm max}}\limits}

\newcommand{\Int}{\displaystyle\int\limits}
\newcommand{\Sup}{\mathop{\rm sup}}

\newcommand{\Lim}{\mathop{{\rm lim}}\limits}

\newcommand{\dist}{\mathop{\rm dist}}
\newcommand{\Img}{\mathop{\rm Im}}
\newcommand{\Real}{\mathop{\rm Re}}

\newcommand{\lal}{{\langle}}
\newcommand{\ral}{{\rangle}}

\newcommand{\be}{\begin{equation}}
\newcommand{\ee}{\end{equation}}
\DeclareMathOperator{\spec}{spec}

\DeclareMathOperator*{\nlim}{{\mathit n}-lim}
\DeclareMathOperator*{\slim}{{\mathit s}-lim}
\DeclareMathOperator*{\wlim}{{\mathit w}-lim}

\newcommand{\ran}{\mathop{\mathrm{Ran}}}
\newcommand{\Ran}{\mathop{\mathrm{Ran}}}
\newcommand{\dom}{\mathop{\mathrm{Dom}}}
\newcommand{\Dom}{\mathop{\mathrm{Dom}}}

\allowdisplaybreaks

\numberwithin{equation}{section}

\newtheorem{theorem}{Theorem}[section]
\newtheorem{corollary}[theorem]{Corollary}
\newtheorem{lemma}[theorem]{Lemma}

\newtheorem{hypothesis}[theorem]{Hypothesis}
\theoremstyle{definition}
\newtheorem{definition}[theorem]{Definition}
\newtheorem{remark}[theorem]{Remark}

\begin{document}
%%%%%%%%%%%%%%%%%%%%%%%%%%%%%%%%%%%%%%%%%%%%%%%%%%%%%%%%%%%%%%%%%

\title[Operator Stieltjes integrals with respect to a spectral measure]
{Operator integrals with respect to a spectral measure and
solutions to some operator equations}

\author[S. Albeverio and A. K. Motovilov]
{Sergio Albeverio  and  Alexander K. Motovilov}

\address{Sergio Albeverio,
Institut f\"ur Angewandte Mathematik, Universit\"at Bonn,
Wegelerstra{\ss}e 6, D-53115 Bonn, Germany; SFB 611, Bonn; BiBoS,
Bielefeld-Bonn; CERFIM, Locarno; Accademia di Architettura, USI, Mendrisio}
\email{albeverio@uni-bonn.de}

\address{Alexander K. Motovilov, Bogoliubov Laboratory of
Theoretical Physics, JINR, Joliot-Curie 6, 141980 Dubna, Moscow
Region, Russia} \email{motovilv@theor.jinr.ru}

%%%\date{{\bf May 08, 2005;}
%%%\textit{Printing date/time}: {\bf \today, \texttime}}

\date{May 08, 2005}



%%%%%%%%%%%%%%%%%%%%%%%%%%%%%%%%%%%%%%%%%%%%%%%%%%%%%%%%
\begin{abstract}
We introduce the concept of Stieltjes integral of an
operator-valued function with respect to the spectral measure
associated with a normal operator. We give sufficient conditions
for the existence of this integral and find bounds on its norm.
The results obtained are applied to the Sylvester and Riccati
operator equations. Assuming that the entry $C$ is a normal
operator, the spectrum of the entry $A$ is separated from the
spectrum of $C$, and $D$ is a bounded operator, we obtain a
representation for the strong solution $X$ to the Sylvester
equation $XA-CX=D$ in the form of an operator Stieltjes integral
with respect to the spectral measure of $C$. By using this result
we then establish sufficient conditions for the existence of a
strong solution to the operator Riccati equation $YA-CY+YBY=D$
where $B$ is another bounded operator.
\end{abstract}
%%%%%%%%%%%%%%%%%%%%%%%%%%%%%%%%%%%%%%%%%%%%%%%%%%%%%

\maketitle

%%%%%%%%%%%%%%%%%%%%%%%%%%%%%%%%%%%%%%%%%%%%%%%%%%%%%
\section{Introduction}
\label{SIntro}
%%%%%%%%%%%%%%%%%%%%%%%%%%%%%%%%%%%%%%%%%%%%%%%%%%%%%%%

In this paper we deal with integrals formally written as
\begin{equation}
\label{FdE}
%J_{\rm l}=
%\cL(F)=
  \int_{\Omega} F(z)dE(z) \text{ \, or \,
}
%\cR(G)=
 \int_{\Omega}dE(z)G(z),
\end{equation}
where $\Omega$ is a Borel subset of the complex plane $\bbC$ and
$E(\cdot)$ the spectral measure associated with a normal operator
on a Hilbert space $\fK$. It is assumed that $F$ is an
operator-valued function on $\Omega$ with values in the space
$\cB(\fK,\fH)$ of bounded operators from $\fK$ to another Hilbert
space $\fH$. Similarly, $G$ is assumed to be an operator-valued
function mapping $\Omega$ into $\cB(\fH,\fK)$. Clearly, a
reasonable definition of integration in \eqref{FdE}
should yield operators from $\fK$ to $\fH$ and from $\fH$ to
$\fK$, respectively.

The integrals of the form \eqref{FdE} are of interest in
itself. But they also arise in many applications, in particular in
the study of spectral subspace perturbation problems (see, e.g.,
\cite{AdLMSr}). Such integrals appear to be a useful tool in
the study of the Sylvester and Riccati operator equations (see
\cite{AdLT,AMM,MM99,MotRem}).

There is an important particular case where $F(z)$ is given by
\begin{equation}
\label{Fphi}
F(z)=\varphi(A,z)\,T
\end{equation}
with $\varphi(\zeta,z)$ a sufficiently nice scalar function of
two complex variables $\zeta$ and $z$, $A$ another normal operator
on $\fH$, and $T$ a bounded operator from $\fK$ to $\fH$. In this
case the integral $\int_{\Omega} F(z)dE(z)$ can be understood as a
double operator Stieltjes integral \cite{BS2003}. More precisely,
\begin{equation}
\int_{\Omega} F(z)dE(z)=\int_{\widetilde{\Omega}}\int_\Omega
\varphi(\zeta,z)d\widetilde{E}(\zeta)TdE(z),
\end{equation}
where $\widetilde{E}(\cdot)$ is the spectral measure associated
with $A$ and $\widetilde{\Omega}=\spec(A)$. Mainly due to the
works by M. Sh. Birman and M. Z. Solomyak \cite{BS67,BS68,BS73}
(see also references cited in \cite{BS2003}) there exists already
a rather comprehensive theory of the double operator Stieltjes
integrals.

However, this is not the case with more general integrals of the
form \eqref{FdE}, the ones that contain operator-valued functions
$F(z)$ and $G(z)$ which cannot be written in terms of
functions of normal operators like in \eqref{Fphi}. To the best of
our knowledge, there is a well established approach to operator
integrals of the type \eqref{FdE}  only in the case where the
spectral measure $E(\cdot)$ is associated with a self-adjoint
operator (see \cite{AdLMSr,AG93,AMM,MM99}, and references
therein). In particular, in \cite{AdLMSr} it was proven that the
integrals \eqref{FdE} make sense as Riemann-Stieltjes integrals
whenever $\Omega$ is a finite interval on the real axis and the
operator-valued functions $F$ and $G$ are Lipschitz on $\Omega$.
Moreover, in this case the integrals \eqref{FdE} exist in the
sense of the uniform operator topology. Some sufficient conditions
for the existence of the improper integrals \eqref{FdE} in the
case where the spectral measure $E(\cdot)$ corresponds to an
unbounded self-adjoint operator are given in \cite{AMM,MM99}.

The present paper is aimed at extending the main concepts and
results of the operator Stieltjes integration theory to the case
where the spectral measure $E(\cdot)$ may correspond to an
arbitrary normal operator. Actually, we consider the integrals
\eqref{FdE} in a somewhat more general setup, admitting the
operator-valued functions $F(\lambda,\mu)$ and $G(\lambda,\mu)$ of
two real variables $\lambda=\Real z$ and $\mu=\Img z$, i.e.
functions that may depend on both $z$ and its conjugate
$\overline{z}$. The integrals \eqref{FdE} are introduced in the
usual way as limits (if they exist) of the corresponding
Riemann-Stieltjes integral sums as the norm of a partition
approaches zero (see Definition \ref{IntDef}). Our main
result is as follows (see Theorem \ref{Integr2}).

Assume that $\Omega$ is a finite rectangle in $\bbC$ and
$F(\lambda,\mu)$ is a $\cB(\fK,\fH)$-valued Lipschitz function
defined for $\lambda+\ri\mu\in\Omega$, that is, there is a
$\gamma_1>0$ such that
\begin{align}
\label{LipIntro}
&\|F(\lambda,\mu)-F(\lambda',\mu')\|\leq \gamma_1\,(|\lambda-\lambda'|+|\mu-\mu'|)\\
\nonumber
&\quad\text{ \,whenever \,}\lambda+\ri\mu\in\Omega\text{\, and \,}
\lambda'+\ri\mu'\in\Omega.
\end{align}
If, in addition, for some $\gamma_2>0$ the function $F$ satisfies the condition
\begin{align}
\label{dLipIntro}
&{\|F(\lambda,\mu)-F(\lambda',\mu)-F(\lambda,\mu')+F(\lambda',\mu')\|
\leq \gamma_2\,|\lambda-\lambda'||\mu-\mu'|}\\
\nonumber
&\qquad\qquad\text{ \,whenever \,}\lambda+\ri\mu\in\Omega\text{\, and \,}
\lambda'+\ri\mu'\in\Omega,
\end{align}
then the operator Stieltjes integral
\begin{equation}
\label{FintI}
\int_\Omega F(\Real z,\Img z)dE(z)
\end{equation}
exists in the sense of the uniform operator topology. The same
statement also holds (see Corollary \ref{Integr21}) for an
operator Stieltjes integral of the form
\begin{equation}
\label{GintI}
\int_\Omega dE(z)G(\Real z,\Img z)
\end{equation}
whenever $G$ is a $\cB(\fH,\fK)$-valued function on $\Omega$
satisfying \eqref{LipIntro} and \eqref{dLipIntro}.

This result allows us to extend the concept of an
operator Stieltjes integral with respect to the spectral measure of
a normal operator $C$ to the case where $\Omega=\spec(C)$, by
taking into account that only the values of $F(\lambda,\mu)$ and
$G(\lambda,\mu)$ for $\lambda+\ri\mu\in\spec(C)$ contribute to
\eqref{FintI} and \eqref{GintI}, respectively. By using this
concept we then obtain an integral representation for the solution
$X$ to the operator Sylvester equation
\begin{equation}
\label{SylInt}
XA-CX=D,
\end{equation}
where $A$ is a closed densely defined possibly
unbounded operator on a Hilbert space $\fH$, $C$ a possibly
unbounded normal operator on a Hilbert space $\fK$, and
$D\in\cB(\fH,\fK)$. In particular, under the assumption that
\begin{equation}
\label{distInt}
\dist\bigl(\spec(A),\spec(C)\bigr)>0
\end{equation}
we prove that $X\in\cB(\fH,\fK)$ is a unique strong solution
to \eqref{SylInt} if and only if it can be represented in
the form of the operator Stieltjes integral
\begin{equation}
\label{EqXInt}
X=\int_{\spec(C)} dE_{C}(z)
D(A-z)^{-1},
\end{equation}
which converges in the sense of the strong operator topology
(see Theorem \ref{SylUnb}). So far, such a representation
was only proven in the case where $C$ is a self-adjoint operator
(see \cite[Theorem 2.14]{AMM}).

We apply the results obtained also to the operator Riccati
equation
\begin{equation}
\label{RicInt}
XA-CX+XBX=D,
\end{equation}
where $B\in\cB(\fK,\fH)$ and the entries $A$, $C$, and $D$ satisfy
the same assumptions as in \eqref{SylInt}. If $X\in\cB(\fH,\fK)$
is a strong solution to \eqref{RicInt}
such that $\spec\bigl(\spec(A+BX),\spec(C)\bigr)>0$, by using
\eqref{EqXInt} one writes the Riccati equation in the equivalent
integral form
\begin{equation}
\label{EqRInt}
X=\int_{\spec(C)} dE_{C}(z)
D(A+BX-z)^{-1}.
\end{equation}
Under the assumption \eqref{distInt}  and additional ``smallness''
assumptions upon $B$ and $D$ we prove the existence of a solution
$X\in\cB(\fH,\fK)$ to the integral equation \eqref{EqRInt}. This
solution to \eqref{EqRInt} also solves the Riccati equation
\eqref{RicInt} (see Theorem \ref{QsolvN}).

The plan of the paper is as follows.

In Section \ref{SInt} we introduce the concept of a Stieltjes
integral of an operator-valued function with respect the spectral
measure of a normal operator and prove the main result (Theorem
\ref{Integr2}) concerning sufficient conditions for the existence
of such integrals. We also derive a norm estimate (Lemma \ref{Lnest})
for these integrals.

In Section \ref{SEnorm} we extend the concept \cite{AMM,MM99} of
the norm of a bounded operator with respect to a spectral measure
to the case where this measure is associated with a normal operator.

In Section \ref{SSyl} we discuss the Sylvester equation
\eqref{SylInt}. In particular, we prove that any weak solution to
this equation is also a strong solution. This result allows us
to present refined versions of the representation theorems
\cite[Section 2]{AMM} for the solution $X$ of the operator Silvester
equation, extending the integral representations for $X$ to the
case where the entry $C$ in \eqref{SylInt} is a normal operator.

Finally, in Section \ref{SecRic} we prove the above mentioned
existence result (Theorem \ref{QsolvN}) for the Riccati equation
\eqref{RicInt}.

We conclude the introduction with the description of some more
notations that are used thro\-ughout the paper. The identity
operator on any Hilbert space $\fK$ is denoted by $I$. If $T$ is a
closed operator on $\fK$, by $\spec(T)$ we always denote the
spectrum of $T$. We will also use the notation
$$
\sigma(T)=\left\{(\lambda,\mu)\in\bbR^2\,|\,\lambda+\ri\mu\in\spec(T)\right\}
$$
for the natural imbedding of $\spec(T)$ into the real plane
$\bbR^2$. The set
\begin{equation*}
\cW(T)=\{z\in\bbC\,|\,\,z=\lal Tx,x\ral\text{ \,
for some\,} x\in\Dom(T), \, \|x\|=1\}
\end{equation*}
is called the numerical range of the operator $T$.

\section{Integral of an operator-valued function with respect
to the spectral measure of a normal operator}
\label{SInt}

Let $\mathcal{A}_\mathrm{Borel}(\bbC)$ denote the algebra of Borel
subsets of the complex plane $\bbC$ and let
$\{E(\Omega)\}_{\Omega\in\mathcal{A}_\mathrm{Borel}(\bbC)}$ be the
spectral family of a (possibly unbounded) normal operator $C$.
Recall that $E(\Omega)$'s are orthogonal projections
in $\fK$ that possess the properties (see, e.g., \cite[\S\,5.1]{BSbook})
\begin{align}
\label{Ecap}
&E(\Omega'\cap\Omega'')=E(\Omega')E(\Omega'') \text{\, for any \,}
\Omega',\Omega''\in\mathcal{A}_\mathrm{Borel}(\bbC),\\
\label{Ecup}
&E(\Omega'\cup\Omega'')=E(\Omega')+E(\Omega'') \text{\, whenever \,}
\Omega',\Omega''\in\mathcal{A}_\mathrm{Borel}(\bbC) \text{\, and \,}
\Omega'\cap\Omega''=\emptyset,
\end{align}
and
\begin{equation}
\label{EmI}
E(\bbC)=E\bigl(\spec(C)\bigr)=I.
\end{equation}
With the spectral family
$\{E(\Omega)\}_{\Omega\in\mathcal{A}_\mathrm{Borel}(\bbC)}$ by
$$
\sE(\Lambda)=E\bigl(\{z=\lambda+\ri\mu\,|\,
(\lambda,\mu)\in\Lambda\subset\bbR^2\}\bigr)
$$
we associate the projection-valued measure
\begin{equation*}
\{\sE(\Lambda)\}_{\Lambda\in\mathcal{A}_\mathrm{Borel}(\bbR^2)}
\end{equation*}
on the algebra $\mathcal{A}_\mathrm{Borel}(\bbR^2)$ of Borel
subsets of $\bbR^2$. Clearly, \eqref{Ecap}, \eqref{Ecup}, and \eqref{EmI} are
equivalent to
\begin{align}
\nonumber
&\sE(\Lambda'\cap\Lambda'')=\sE(\Lambda')\sE(\Lambda'') \text{\, for any \,}
\Lambda',\Lambda''\in\mathcal{A}_\mathrm{Borel}(\bbR^2),\\
\nonumber
&\sE(\Lambda'\cup\Lambda'')=\sE(\Lambda')+\sE(\Lambda'') \text{\, whenever \,}
\Lambda',\Lambda''\in\mathcal{A}_\mathrm{Borel}(\bbR^2) \text{\, and \,}
\Lambda'\cap\Lambda''=\emptyset,
\end{align}
and
\begin{equation*}
E(\bbR^2)=E\bigl(\sigma(C)\bigr)=I,
\end{equation*}
respectively.

In terms of the spectral measures $E(\cdot)$ and $\sE(\cdot)$, we write the
spectral decomposition of the normal operator $C$ either as
\begin{equation*}
C=\int_{\bbC}z \, dE(z)=\int_{\spec(C)}z \, dE(z)
\end{equation*}
or as
\begin{equation*}
C=\int_{\bbR^2}(\lambda+\ri\mu)\, d\sE(\lambda,\mu)=
\int_{\sigma(C)}(\lambda+\ri\mu)\, d\sE(\lambda,\mu).
\end{equation*}

Further, we introduce the projection-valued function
$\hsE(\lambda,\mu)$ on $\bbR^2$ by
\begin{equation*}
\hsE(\lambda,\mu)=\sE(\{(x,y)\in\bbR^2\,|\,x<\lambda,\,y<\mu\}).
\end{equation*}
In the following the function $\hsE(\lambda,\mu)$ is called the
spectral function of the normal operator $C$. In contrast to the
case of self-adjoint operators the spectral function of a normal
operator is a function of two real variables.

Clearly, for $\lambda\leq\lambda'$ and $\mu\leq\mu'$ the
additivity property of a spectral measure implies
$$
\hsE(\lambda',\mu')=\hsE(\lambda,\mu)
+\sE\bigl((-\infty,\lambda)\times[\mu,\mu')\bigr)
+\sE\bigl([\lambda,\lambda')\times(-\infty,\mu)\bigr)
+\sE\bigl([\lambda,\lambda')\times[\mu,\mu')\bigr)
$$
and hence
\begin{equation}
\label{Enon} \hsE(\lambda,\mu)\leq \hsE(\lambda',\mu')\text{ if }
\lambda\leq\lambda'\text{ and }\mu\leq\mu'\,,
\end{equation}
that is, $\hsE(\lambda,\mu)$ is a non-decreasing function in both
variables $\lambda$ and $\mu$.

One also observes that if $\lambda\leq\lambda'$ and $\mu\leq\mu'$
then
\begin{equation}
\label{dpos1}
\hsE(\lambda',\mu')-\hsE(\lambda',\mu)-\hsE(\lambda,\mu')+
\hsE(\lambda,\mu)=E\bigl([\lambda,\lambda')\times[\mu,\mu')\bigr)
\end{equation}
and thus
\begin{equation}
\label{dpos2}
\begin{array}{c}
 \hsE(\lambda',\mu)-\hsE(\lambda,\mu) \leq
\hsE(\lambda',\mu')-\hsE(\lambda,\mu') \\
\text{ for any } \mu,\mu'\in\bbR\text{ such that } \mu\leq\mu'
\text{ and any } \lambda,\lambda'\in\bbR \text{ such that }
\lambda\leq\lambda'.
\end{array}
\end{equation}

We remark that \eqref{Enon} implies
\begin{equation}
\label{EnonX} \|\hsE(\lambda',\mu')x\|\leq
\|\hsE(\lambda,\mu)x\|\text{ for any } x\in\fK \,\text{ whenever  }
\lambda'\leq\lambda\text{ and } \mu'\leq\mu\,.
\end{equation}
Indeed
$$
\|\hsE(\lambda,\mu)x\|^2=\langle\hsE(\lambda,\mu)x,\hsE(\lambda,\mu)x\rangle=
\langle\hsE(\lambda,\mu)x,x\rangle
$$
and
$$
\|\hsE(\lambda',\mu')x\|^2=\langle\hsE(\lambda',\mu')x,\hsE(\lambda',\mu')x\rangle=
\langle\hsE(\lambda',\mu')x,x\rangle.
$$
Hence \eqref{EnonX} is a consequence of \eqref{Enon}.

By using the fact that both right-hand and left-hand sides of
\eqref{dpos2} are orthogonal projections, in a similar way one
concludes that
\begin{equation}
\label{dposX}
\begin{array}{c}
\left\|\bigl(\hsE(\lambda',\mu)-\hsE(\lambda,\mu)\bigr)x\right\|\leq
\left\|\bigl(\hsE(\lambda',\mu')-\hsE(\lambda,\mu')\bigr)x\right\|\\
\text{ for any } x\in\fK \text{ whenever } \mu\leq\mu' \text{ and
}  \lambda\leq\lambda'.
\end{array}
\end{equation}

For notational setup we adopt the following

\begin{hypothesis}
\label{Hypo1} Let $\fH$ and $\fK$ be Hilbert spaces,
$\Delta=[a,b)\times[c,d)$ a rectangle in $\bbR^2$ with
$-\infty<a<b<+\infty$ and $-\infty<c<d<+\infty$. Let
$\{\sE(\Lambda)\}_{\Lambda\in\mathcal{A}_\mathrm{Borel}(\bbR^2)}$
be the spectral family associated with a normal operator on $\fK$.
\end{hypothesis}

%%%%%%%%%%%%%%%%%%%%%%%%%%%%%%%%%%%%%%%%%%%%%%%%%%%%%%%%%%%%%%%%%
\begin{definition}
\label{IntDef} Assume Hypothesis \ref{Hypo1}.
An operator-valued function
$$
F:\,\Delta\to\cB(\fK,\fH)
$$
is said to be uniformly {\rm(}resp. strongly, weakly{\rm)}
integrable from the right with respect the spectral measure
$d\sE(\lambda,\mu)$ on $\Delta$ if the limit
\begin{equation}
\label{defIntR}
\Int_\Delta F(\lambda,\mu)\,d\sE(\lambda,\mu)=
\Lim_{\Max_{j=1}^m |\delta_j^{(m)}|+\Max_{k=1}^n
|\omega_k^{(n)}|\rightarrow0}\,\,
\sum\limits_{j=1}^m \sum\limits_{k=1}^n  F(\xi_j,\zeta_k)\,
\sE(\delta_j^{(m)}\times\omega_k^{(n)})
\end{equation}
exists in the uniform {\rm(}resp. strong, weak{\rm)} operator
topology. Here, $\delta_j^{(m)}=[\lambda_{j-1},\lambda_j)$,
$j=1,2,\ldots,m$,  and $\omega_k^{(n)}=[\mu_{k-1},\mu_k)$,
$k=1,2,\ldots,n$, where $a=\lambda_0<\lambda_1<\ldots<\lambda_m=b$
and $c=\mu_0<\mu_1<\ldots<\mu_n=b$ are partitions of the intervals
$[a,b)$ and $[c,d)$, respectively;  $\xi_j\in \delta_j^{(m)}$ and
$\zeta_k\in \delta_k^{(n)}$ are arbitrarily chosen points,
$|\delta_k^{(j)}|=\lambda_j-\lambda_{j-1}$ and
$|\omega_k^{(n)}|=\mu_k-\mu_{k-1}$. The limit
value~{\rm(\Ref{defIntR})}, if it exists, is called the right
Stieltjes integral of the operator-valued function $F$ with respect to the
measure $d\sE(\lambda,\mu)$ on $\Delta$.

Similarly, a function
$$
G:\,\Delta\rightarrow\cB(\fH,\fK)
$$
is said to be uniformly {\rm(}resp. strongly, weakly{\rm)}
integrable from the left with respect to the measure $d\sE(\lambda,\mu)$ on
$\Delta$, if the limit
\begin{equation}
\label{defIntL}
\Int_\Delta d\sE(\lambda,\mu)\,G(\lambda,\mu)=
\Lim_{\Max_{j=1}^m |\delta_j^{(m)}|+\Max_{k=1}^n
|\omega_k^{(n)}|\rightarrow0}\,\,
\sum\limits_{j=1}^m \sum\limits_{k=1}^n
\sE(\delta_j^{(m)}\times\omega_k^{(n)})\,G(\xi_j,\zeta_k)\,
\end{equation}
exists in the uniform {\rm(}resp. strong, weak{\rm)} operator
topology. The corresponding limit value {\rm(\Ref{defIntL})}, if
it exists, is called the left Stieltjes integral of the
operator-valued function $G$ with respect to the measure
$d\sE(\lambda,\mu)$ on $\Delta$.
\end{definition}

The following statement can be considered an extension of
\cite[Lemma 10.5]{MM99} to the case of the spectral measure
associated with a normal operator.

%%%%%%%%%%%%%%%%%%%%%%%%%%%%%%%%%%%%%%%%%%%%%%%%%%%%%%%%%%%%%%%%%%
\begin{lemma}
\label{Integr1} Assume Hypothesis \ref{Hypo1}. Then an operator-valued
function $F(\lambda,\mu)$,
$$
F:\,\Delta\rightarrow\cB(\fK,\fH),
$$
is integrable in the weak {\rm(}uniform{\rm)} operator
topology with respect to the measure $d\sE(\lambda,\mu)$ on the rectangle
$\Delta$ from the left if and only if  the function
$[F(\lambda,\mu)]^*$ is integrable in the weak
{\rm(}uniform{\rm)} operator topology with respect to the measure
$d\sE(\lambda,\mu)$ on $\Delta$ from the right and then
\begin{equation}
\label{JJadj}
\left[\int_\Delta F(\lambda,\mu)\,d\sE(\lambda,\mu)\right]^*
     =\int_\Delta d\sE(\lambda,\mu)\,[F(\lambda,\mu)]^*.
\end{equation}
\end{lemma}
%%%%%%%%%%%%%%%%%%%%%%%%%%%%%%%%%%%%%%%%%%%%%%%%%%%%%%%%%%%%%%
\begin{proof}
Like in the proof of Lemma 10.5 in \cite{MM99}, the assertion is
proven by taking into account the continuity property of the
involution $T\rightarrow T^*$ with respect to operator uniform and
operator weak convergence in $\cB(\fK,\fH)$. It suffices to apply this
property to the integral sums in \eqref{defIntR} and
\eqref{defIntL}.
\end{proof}
\begin{remark}
\label{RemConv} Since the involution $T\to T^*$ is not continuous
with respect to the strong convergence (see, e.g., \cite[\S 2.5]{BSbook}), the convergence of one of the integrals
\eqref{JJadj} in the strong operator topology in general only
implies the convergence of the other one in the weak operator
topology.
\end{remark}
%%%%%%%%%%%%%%%%%%%%%%%%%%%%%%%%%%%%%%%%%%%%%%%%%%%%%%%%%%%%%%

Some sufficient conditions for the integrability of an
operator-valued function $F(\lambda,\mu)$ with respect to the spectral
measure of a normal operator are described in the following
statement.

%%%%%%%%%%%%%%%%%%%%%%%%%%%%%%%%%%%%%%%%%%%%%%%%%%%%%%%%%%%%%%
\begin{theorem}
\label{Integr2} Assume Hypothesis \ref{Hypo1}. Suppose that for
some $\gamma_1>0$ the operator-valued function
$F:\,\Delta\rightarrow\cB(\fK,\fH)$ satisfies the Lipschitz
condition
\begin{eqnarray}
\label{Lipschitz}
{\|F(\lambda,\mu)-F(\lambda',\mu')\|\leq \gamma_1\,(|\lambda-\lambda'|+|\mu-\mu'|)}
\end{eqnarray}
and for some $\gamma_2>0$ the condition
\begin{equation}
\label{dLipschitz}
{\|F(\lambda,\mu)-F(\lambda',\mu)-F(\lambda,\mu')+F(\lambda',\mu')\|
\leq \gamma_2\,|\lambda-\lambda'||\mu-\mu'|}
\end{equation}
$$
\text{for any}\quad \lambda,\lambda'\in[a,b)\quad \text{and}\quad \mu,\mu'\in[c,d).
$$
Then the function $F$ is right-integrable on $\Delta$ with respect
to the spectral measure $d\sE(\lambda,\mu)$ in the sense of the
uniform operator topology.
\end{theorem}
%%%%%%%%%%%%%%%%%%%%%%%%%%%%%%%%%%%%%%%%%%%%%%%%%%%%%%%%%%%%%%%%
\begin{proof}
Let $\{\delta_j^{(m)}\}_{j=1}^{m}$ and
$\{\omega_k^{(n)}\}_{k=1}^n$ be partitions of the intervals
$[a,b)$ and $[c,d)$, respectively, and let
$\Delta^{(m,n)}_{jk}=\delta_k^{(m)}\times\omega_k^{(n)}$,
$j=1,2,\ldots,m,\, k=1,2,\ldots,n$. Assume that
$\{(\xi_j,\zeta_k)\in\Delta^{(m,n)}_{jk},\,\, j=1,2,\ldots,m,\,\,
k=1,2,\ldots,n\}$ and
$\{(\xi'_j,\zeta'_k)\in\Delta^{(m,n)}_{jk}, \,\,
j=1,2,\ldots,m,\,\, k=1,2,\ldots,n\}$ are two sets of points.

First, we prove that the limit \eqref{defIntR} (if it exists) does not depend on
the choice of the points $(\xi_j,\zeta_k)$ within the partition rectangles
$\Delta^{(m,n)}_{jk}$.

Let
\begin{equation}
\label{Jmn} J_{m,n}=\sum_{j=1}^m \sum_{k=1}^n F(\xi_j,\zeta_k)
\sE(\Delta^{(m,n)}_{jk})
\end{equation}
and
$$
J'_{m,n}=\sum_{j=1}^m \sum_{k=1}^n
F(\xi'_j,\zeta'_k)\sE(\Delta^{(m,n)}_{jk}).
$$
Observe that
\begin{equation}
\label{EE4}
\sE(\Delta^{(m,n)}_{jk})=\hsE(\lambda_j,\mu_k)-\hsE(\lambda_{j-1},\mu_k)
-\hsE(\lambda_{j},\mu_{k-1})+\hsE(\lambda_{j-1},\mu_{k-1})
\end{equation}
and hence
\begin{align*}
J_{m,n}-J'_{m,n}=&\sum_{j=1}^m \sum_{k=1}^n
\bigl([F(\xi_j,\zeta_k)-F(\xi'_j,\zeta'_k)\bigr]\,\\
&\qquad\times\bigl[\hsE(\lambda_j,\mu_k)-\hsE(\lambda_{j-1},\mu_k)
-\hsE(\lambda_{j},\mu_{k-1})+\hsE(\lambda_{j-1},\mu_{k-1})\bigr].
\end{align*}
Represent the difference $J_{m,n}-J'_{m,n}$ as  the sum of two
terms that are more convenient for estimating
\begin{equation}
\label{JJLL}
J_{m,n}-J'_{m,n}=L_1+L_2,
\end{equation}
where
\begin{align*}
L_1=&\sum_{j=1}^m \sum_{k=1}^n
\bigl([F(\xi_j,\zeta_k)-F(\xi'_j,\zeta_k)\bigr]\,\\
&\qquad\times\bigl[\hsE(\lambda_j,\mu_k)-\hsE(\lambda_{j-1},\mu_k)
-\hsE(\lambda_{j},\mu_{k-1})+\hsE(\lambda_{j-1},\mu_{k-1})\bigr]
\end{align*}
and
\begin{align}
\label{JJL2}
L_2=&\sum_{j=1}^m \sum_{k=1}^n
\bigl([F(\xi'_j,\zeta_k)-F(\xi'_j,\zeta'_k)\bigr]\,\\
\nonumber
&\qquad\times\bigl[\hsE(\lambda_j,\mu_k)-\hsE(\lambda_{j-1},\mu_k)
-\hsE(\lambda_{j},\mu_{k-1})+\hsE(\lambda_{j-1},\mu_{k-1})\bigr].
\end{align}

By inspection one verifies that
\begin{align}
\label{J-J}
&L_1=S_1+S_2+S_3,
\end{align}
where
\begin{align}
\nonumber
S_1&=\sum_{j=1}^m \bigl[F(\xi'_j,\zeta_1)-F(\xi_j,\zeta_1)\bigr]
\bigl[\hsE(\lambda_j,\mu_0)-\hsE(\lambda_{j-1},\mu_0)\bigr],\\
\nonumber
S_2&=\sum_{j=1}^m \bigl[F(\xi_j,\zeta_n)-F(\xi'_j,\zeta_n)\bigr]
\bigl[\hsE(\lambda_j,\mu_n)-\hsE(\lambda_{j-1},\mu_n)\bigr],
\end{align}
and
\begin{align*}
S_3=&\sum_{j=1}^m \,\,\sum_{k=1}^{n-1}
\bigl[F(\xi_j,\zeta_k)-F(\xi'_j,\zeta_k)-
F(\xi_j,\zeta_{k+1})+F(\xi'_j,\zeta_{k+1})\bigr]\\
&\qquad\qquad\times
\bigl[\hsE(\lambda_j,\mu_k)-\hsE(\lambda_{j-1},\mu_k)\bigr].
\end{align*}
Clearly, for any $x\in\fK$ by the Lipschitz property \eqref{Lipschitz}
$$
\|S_1 x\|\leq \gamma_1 \sum_{j=1}^m |\delta^{(m)}_j|
\|\bigl(\hsE(\lambda_j,\mu_0)-\hsE(\lambda_{j-1},\mu_0)\bigr)x\|
$$
and hence
\begin{align}
\nonumber
\|S_1 x\|&\leq \gamma_1\left(\sum_{j=1}^m |\delta^{(m)}_j|^2\right)^{1/2}
\left(\sum_{j=1}^m \lal\bigl(
\hsE(\lambda_j,\mu_0)-\hsE(\lambda_{j-1},\mu_0)\bigr)x,x \ral\right)^{1/2}\\
\nonumber
&\leq \gamma_1 \left(\Max_{j=1}^m|\delta^{(m)}_j|\right)^{1/2} \sqrt{b-a}
\left(\lal\bigl(\hsE(b,\mu_0)-\hsE(a,\mu_0)\bigr)x,x\ral\right)^{1/2}\\
\label{S1}
&\leq \gamma_1 \left(\Max_{j=1}^m|\delta^{(m)}_j|\right)^{1/2} \sqrt{b-a}
\left\|\hsE(b,c)-\hsE(a,c)\right\|\|x\|,
\end{align}
since
\begin{align}
\nonumber
\lal\bigl(\hsE(b,\mu_0)-\hsE(a,\mu_0)\bigr)x,x\ral
&=\lal\bigl(\hsE(b,\mu_0)-\hsE(a,\mu_0)\bigr)^2x,x\ral \\
\nonumber
&=\lal\bigl(\hsE(b,\mu_0)-\hsE(a,\mu_0)\bigr)x,\bigl(\hsE(b,\mu_0)-\hsE(a,\mu_0)\bigr)x\ral\\
\label{nxest}
&=\left\|\bigl(\hsE(b,\mu_0)-\hsE(a,\mu_0)\bigr)x\right\|^2
\end{align}
and $\mu_0=c$.

In a similar way one shows that for any $x\in\fK$
\begin{equation}
\label{S2}
\|S_2x\|\leq \gamma_1\left(\Max_{j=1}^m|\delta^{(m)}_j|\right)^{1/2}\sqrt{b-a}
\left\|\hsE(b,d)-\hsE(a,d)\right\|\|x\|.
\end{equation}

Finally, by using \eqref{dLipschitz} at the fist step,  for
any $x\in\fK$ one estimates $S_3x$ as follows:
\begin{align}
\nonumber
\|S_3x\|&\leq \gamma_2\sum_{j=1}^m\sum_{k=1}^{n-1}|\delta^{(m)}_j||\omega^{(n)}_k|
\left\|\bigl(\hsE(\lambda_j,\mu_k)-\hsE(\lambda_{j-1},\mu_k)\bigr)x\right\|\\
\nonumber
&\qquad = \gamma_2\sum_{k=1}^{n-1} |\omega^{(n)}_k| \sum_{j=1}^m |\delta^{(m)}_j|
\left\|\bigl(\hsE(\lambda_j,\mu_k)-\hsE(\lambda_{j-1},\mu_k)\bigr)x\right\|\\
\nonumber
& \leq\gamma_2\sum_{k=1}^{n-1} |\omega^{(n)}_k| \left(\sum_{j=1}^m |\delta^{(m)}_j|^2\right)^{1/2}
\left(\sum_{j=1}^m \lal\bigl(\hsE(\lambda_j,\mu_k)-\hsE(\lambda_{j-1},\mu_k)\bigr)x,x\ral
\right)^{1/2}\\
\nonumber
&\leq \gamma_2\,(d-c)\left(\Max_{j=1}^m|\delta^{(m)}_j|\right)^{1/2}\sqrt{b-a}\,\,\,
\Max_{k=1}^{n-1}\left(\lal(\hsE(b,\mu_k)-\hsE(a,\mu_k)\bigr)x,x\ral\right)^{1/2}\\
\nonumber
&\leq \gamma_2\,\sqrt{b-a}\,(d-c)\,\,\left(\Max_{j=1}^m|\delta^{(m)}_j|\right)^{1/2}\,\Max_{k=1}^{n-1}
\|\hsE(b,\mu_k)-\hsE(a,\mu_k)\|\, \|x\|
\end{align}
by applying \eqref{nxest} (with $\mu_0$ replaced by $\mu_k$) at
the last step. Obviously, by \eqref{dpos2} for any $k=1,2,\ldots,n-1$
$$
\hsE(b,\mu_k)-\hsE(a,\mu_k)\leq \hsE(b,\mu_n)-\hsE(a,\mu_n)=\hsE(b,d)-\hsE(a,d)
$$
and we arrive at the following final estimate for $S_3$:
\begin{align}
\label{S3}
\|S_3 x\|\leq \gamma_2\,\sqrt{b-a}\,(d-c)\,\left(\Max_{j=1}^m |\delta^{(m)}_j|\right)^{1/2}
\,
\|\hsE(b,d)-\hsE(a,d)\|\, \|x\|.
\end{align}

Combining \eqref{J-J}, \eqref{S1}, \eqref{S2}, and \eqref{S3}
one concludes that for any $x\in\fK$
\begin{align}
\label{L1est}
\|L_1x\|\leq M_1 \left(\Max_{j=1}^m |\delta^{(m)}_j|\right)^{1/2}\,\sqrt{b-a}\,\|x\|,
\end{align}
where
$$
M_1=\gamma_1\|\hsE(b,c)-\hsE(a,c)\|+\bigl(\gamma_1+\gamma_2(d-c)\bigr)\|
\hsE(b,d)-\hsE(a,d)\|.
$$

In a similar way one proves that an analogous estimate holds for the term
$L_2$ given by \eqref{JJL2},
\begin{align}
\label{L2est}
\|L_2x\|\leq  M_2 \left(\Max_{k=1}^n |\omega^{(n)}_k|\right)^{1/2}\,
\sqrt{d-c}\,\|x\|,
\end{align}
where
$$
M_2=\gamma_1\|\hsE(a,d)-\hsE(a,c)\|+\bigl(\gamma_1+\gamma_2(b-a)\bigr)\|
\hsE(b,d)-\hsE(b,c)\|.
$$

Combining \eqref{JJLL}, \eqref{L1est}, and \eqref{L2est} proves
that if the sum in \eqref{defIntR} converges strongly (resp.
weakly, in the operator norm topology) for some choice of the
numbers $\{\xi_j\in\delta^{(m)}_j\}_{j=1}^m$ and
$\{\zeta_k\in\omega^{(n)}_k\}_{k=1}^n$, then it converges strongly
(resp. weakly, in the operator norm topology) to the same limit
for any other choice of these numbers, in particular this takes
place for the choice where
\begin{equation}
\label{xiz} \xi_j=\lambda_{j-1}, \quad j=1,2,\ldots,m, \quad
\text{and} \quad \zeta_k=\mu_{j-1}, \quad k=1,2,\ldots,n.
\end{equation}
It remains to prove that the double sequence of the operators
$J_{m,n}$ given by \eqref{Jmn} has a limit as $m\to\infty$,
$n\to\infty$, $\Max_{j=1}^m|\delta^{(m)}_j|\to0$, and
$\Max_{k=1}^n|\omega^{(n)}_k|\to0$.

Assume that there are two different partitions of the interval
$[a,b)$ containing $m$ and $\widetilde{m}>m$ subintervals, respectively,
and two different partitions of the interval $[c,d)$ containing
$n$ and $\widetilde{n}>n$ subintervals, respectively Further, assume, without
loss of generality, that all the points of the $m$-partition of the
interval $[a,b)$ are the points of the $\widetilde{m}$-partition
of $[a,b)$ and all the points of the $\widetilde{n}$-partition of
the interval $[c,d)$ are the points of the $n$-partition of
$[c,d)$. Denote by $\lambda_{j,s}$, $s=0,1,\ldots,p_j,$ the points
of the $\widetilde{m}$-partition of $[a,b)$ that belong to the
interval $[\lambda_{j-1},\lambda_j]$ and by $\mu_{k,t}$,
$t=0,1,\ldots,q_k,$ the points of the $\widetilde{n}$-partition of
$[c,d)$ belonging to the interval $[\mu_{k-1},\mu_k]$. By
definition of partition subintervals we have
\begin{align}
\label{lamj}
\lambda_{j-1}=\lambda_{j,0}<\lambda_{j,1}<&\ldots<\lambda_{j,p_j}=\lambda_j,\quad
j=1,2,\ldots,m,\\
\label{muk}
\mu_{k-1}=\mu_{k,0}<\mu_{k,1}<&\ldots<\mu_{k,q_k}=\mu_k,\quad
k=1,2,\ldots,n.
\end{align}

In the following we assume that one chooses the points $\xi_j$,
$\zeta_k$ and $\xi_{j,s}$, $\zeta_{k,t}$ like in \eqref{xiz} and
thus
\begin{equation}
\label{JmnL}
J_{m,n}=\sum_{j=1}^m \sum_{k=1}^n F(\lambda_{j-1},\mu_{k-1}) \sE(\Delta^{(m,n)}_{jk})
\end{equation}
and
\begin{equation}
\label{tJmnL}
J_{\widetilde{m},\widetilde{n}}=
\sum_{j=1}^m \sum_{k=1}^n \sum_{s=1}^{p_j}\sum_{t=1}^{q_k}
F(\lambda_{j,s-1},\mu_{k,t-1}) \sE(\Delta^{(\widetilde{m},\widetilde{n})}_{j,s;k,t}),
\end{equation}
where
$$
\Delta^{(\widetilde{m},\widetilde{n})}_{j,s;k,t}=[\lambda_{j,s-1},\lambda_{j,s})
\times[\mu_{k,t-1;k,t}).
$$
By \eqref{EE4}
\begin{align}
\nonumber
J_{m,n}&=\sum_{j=1}^m \sum_{k=1}^n \,F(\lambda_{j-1},\mu_{k-1})\\
\label{JmnLe}
&\quad\times\left[\hsE(\lambda_j,\mu_k)-\hsE(\lambda_{j-1},\mu_k)
-\hsE(\lambda_{j},\mu_{k-1})+\hsE(\lambda_{j-1},\mu_{k-1})\right]
\end{align}
and
\begin{align}
\nonumber
J_{\widetilde{m},\widetilde{n}}&=
\sum_{j=1}^m \sum_{k=1}^n \sum_{s=1}^{p_j}\sum_{t=1}^{q_k}
F(\lambda_{j,s-1},\mu_{k,t-1})\\
\label{tJmnLe}
&\quad\times\left[\hsE(\lambda_{j,s},\mu_{k,t})-\hsE(\lambda_{j,s-1},\mu_{k,t})
-\hsE(\lambda_{j,s},\mu_{k,t-1})+\hsE(\lambda_{j,s-1},\mu_{k,t-1})\right].
\end{align}
Taking into account \eqref{lamj} and \eqref{muk} one verifies by
inspection that
\begin{equation}
\label{JJdif}
J_{\widetilde{m},\widetilde{n}}-J_{mn}=\sum_{j=1}^m
\sum_{k=1}^n (T^{(1)}_{jk}+T^{(2)}_{jk}+T^{(3)}_{jk}),
\end{equation}
where
\begin{align}
\nonumber
T^{(1)}_{jk}&=\sum_{s=1}^{p_j-1}
\left[F(\lambda_{j,s},\mu_{k,0})-F(\lambda_{j,s-1},\mu_{k,0})\right]\\
\nonumber
&\qquad\times\left[\hsE(\lambda_{j,p_j},\mu_{k,q_k})-\hsE(\lambda_{j,s},\mu_{k,q_k})
-\hsE(\lambda_{j,p_j},\mu_{k,0})+\hsE(\lambda_{j,s},\mu_{k,0})\right]\\
\nonumber &=\sum_{s=1}^{p_j-1}
\left[F(\lambda_{j,s},\mu_{k-1})-F(\lambda_{j,s-1},\mu_{k-1})\right]\\
\nonumber
&\qquad\times\left[\hsE(\lambda_{j},\mu_{k})-\hsE(\lambda_{j,s},\mu_{k})
-\hsE(\lambda_{j},\mu_{k-1})+\hsE(\lambda_{j,s},\mu_{k-1})\right],
\end{align}
\begin{align}
\nonumber T^{(2)}_{jk}&=\sum_{t=1}^{q_k-1}
\left[F(\lambda_{j,0},\mu_{k,t})-F(\lambda_{j,0},\mu_{k,t-1})\right]\\
\nonumber
&\qquad\times\left[\hsE(\lambda_{j,p_j},\mu_{k,q_k})-\hsE(\lambda_{j,0},\mu_{k,q_k})
-\hsE(\lambda_{j,p_j},\mu_{k,t})+\hsE(\lambda_{j,0},\mu_{k,t})\right]\\
\nonumber &=\sum_{t=1}^{q_k-1}
\left[F(\lambda_{j-1},\mu_{k,t})-F(\lambda_{j-1},\mu_{k,t-1})\right]\\
\nonumber
&\qquad\times\left[\hsE(\lambda_{j},\mu_k)-\hsE(\lambda_{j-1},\mu_{k})
-\hsE(\lambda_{j},\mu_{k,t})+\hsE(\lambda_{j-1},\mu_{k,t})\right],
\end{align}
and
\begin{align}
\nonumber
T^{(3)}_{jk}&=\sum_{s=1}^{p_j-1}\,\sum_{t=1}^{q_k-1}
\left[F(\lambda_{j,s},\mu_{k,t})-F(\lambda_{j,s-1},\mu_{k,t})-
F(\lambda_{j,s},\mu_{k,t-1})+F(\lambda_{j,s-1},\mu_{k,t-1})\right]\\
\label{T3jk}
&\qquad\times\left[\hsE(\lambda_{j},\mu_{k})-\hsE(\lambda_{j,s},\mu_{k})
-\hsE(\lambda_{j},\mu_{k,t})+\hsE(\lambda_{j,s},\mu_{k,t})\right].
\end{align}
One also notes that
\begin{align}
\nonumber
\sum_{k=1}^n T^{(1)}_{jk}=&\sum_{s=1}^{p_j-1}\left\{\left[F(\lambda_{j,s},\mu_{n-1})-F(\lambda_{j,s-1},\mu_{n-1})\right]
\left[\hsE(\lambda_j,\mu_n)-\hsE(\lambda_{j,s},\mu_n)\right]\right.\\
\nonumber
&\qquad-\left[F(\lambda_{j,s},\mu_0)-F(\lambda_{j,s-1},\mu_0)\right]
\left[\hsE(\lambda_j,\mu_0)-\hsE(\lambda_{j,s},\mu_0)\right]\\
\nonumber
&\qquad+\sum_{k=1}^{n-1}\left[F(\lambda_{j,s},\mu_{k-1})-F(\lambda_{j,s-1},\mu_{k-1})
-F(\lambda_{j,s},\mu_k)+F(\lambda_{j,s-1},\mu_k)\right]\\
\label{T1jk}
&\left.\qquad\qquad\times\left[\hsE(\lambda_j,\mu_k)-\hsE(\lambda_{j,s},\mu_k)\right]\right\}
\end{align}
and
\begin{align*}
\sum_{j=1}^m T^{(2)}_{jk}=&\sum_{t=1}^{q_k-1}\left\{\left[F(\lambda_{m-1},\mu_{k,t})-F(\lambda_{m-1},\mu_{k,t-1})\right]
\left[\hsE(\lambda_m,\mu_k)-\hsE(\lambda_m,\mu_{k,t})\right]\right.\\
&\qquad-\left[F(\lambda_0,\mu_{k,t})-F(\lambda_0,\mu_{k,t-1})\right]
\left[\hsE(\lambda_0,\mu_k)-\hsE(\lambda_0,\mu_{k,t})\right]\\
&\qquad+\sum_{j=1}^{m-1}\left[F(\lambda_{j-1},\mu_{k,t})-F(\lambda_{j-1},\mu_{k,t-1})
-F(\lambda_j,\mu_{k,t})+F(\lambda_j,\mu_{k,t-1})\right]\\
&\left.\qquad\qquad\times\left[\hsE(\lambda_j,\mu_k)-\hsE(\lambda_j,\mu_{k,t})\right]\right\}.
\end{align*}

By \eqref{Lipschitz} and \eqref{dLipschitz} the equality
\eqref{T1jk} implies that for any $x\in\fK$
\begin{align}
\nonumber
\biggl\|\sum_{k=1}^n T^{(1)}_{jk}x\biggr\|\leq
\sum_{s=1}^{p_j-1}&\left\{\gamma_1\,|\lambda_{j,s}-\lambda_{j,s-1}|
\left\|\bigl(\hsE(\lambda_j,\mu_n)-\hsE(\lambda_{j,s},\mu_n)\bigr)x\right\|
\phantom{\sum_{s=1}^n}\right.\\
\nonumber
&\quad+\gamma_1\,|\lambda_{j,s}-\lambda_{j,s-1}|
\left\|\bigl(\hsE(\lambda_j,\mu_0)-\hsE(\lambda_{j,s},\mu_0)\bigr)x\right\|\\
\nonumber
&\quad+\left.\sum_{k=1}^{n-1}\gamma_2\,|\lambda_{j,s}-\lambda_{j,s-1}|
\,|\omega_k^{(n)}|\left\|\bigl(\hsE(\lambda_j,\mu_k)-
\hsE(\lambda_{j,s},\mu_k)\bigr)x\right\|\right\}
\end{align}
and thus
\begin{align}
\nonumber \biggl\|\sum_{k=1}^n T^{(1)}_{jk}x\biggr\|
&leq\gamma_1\,|\delta^{(m)}_j|\Max_{s=1}^{p_j-1}
\left\|\bigl(\hsE(\lambda_j,d)-\hsE(\lambda_{j,s},d)\bigr)x\right\|\\
\nonumber &\quad +\gamma_1\,|\delta^{(m)}_j|\Max_{s=1}^{p_j-1}
\left\|\bigl(\hsE(\lambda_j,c)-\hsE(\lambda_{j,s},c)\bigr)x\right\|\\
\label{Tjk1x} & \quad +\gamma_2\,|\delta^{(m)}_j| (d-c)|
\Max_{k=1}^{n-1}\,\Max_{s=1}^{p_j-1}
\left\|\bigl(\hsE(\lambda_j,\mu_k)-
\hsE(\lambda_{j,s},\mu_k)\bigr)x\right\|,
\end{align}
since $\mu_0=c$, $\mu_n=d$,
$$
\sum_{s=1}^{p_j-1}|\lambda_{j,s}-\lambda_{j,s-1}|=
\lambda_{j,p_{j-1}}-\lambda_{j-1}<\lambda_j-\lambda_{j-1}=|\delta^{(m)}_j|,
$$
and $\sum\limits_{k=1}^{n-1}|\omega_k^{(n)}|<d-c$. Obviously
\begin{align}
\label{Edif}
\left\|\bigl(\hsE(\lambda_j,\mu)-\hsE(\lambda_{j,s},\mu)\bigr)x\right\|&\leq
\left\|\bigl(\hsE(\lambda_j,\mu)-
\hsE(\lambda_{j-1},\mu)\bigr)x\right\|\text{ for any }\mu\in\bbR
\end{align}
since
\begin{align}
\nonumber
\left\|\bigl(\hsE(\lambda_j,\mu)-\hsE(\lambda_{j,s},\mu)\bigr)x\right\|^2=&
\bigl\langle
\sE\bigl([\lambda_{j,s},\lambda_j)\times(-\infty,\mu)\bigr)x,x\bigr\rangle\\
\nonumber &\quad \leq  \bigl\langle
\sE\bigl([\lambda_{j-1},\lambda_j)\times(-\infty,\mu)\bigr)x,x\bigr\rangle\\
\nonumber
&\quad\qquad=\left\|\bigl(\hsE(\lambda_j,\mu)-\hsE(\lambda_{j-1},\mu)
\bigr)x\right\|^2.
\end{align}
Then by using \eqref{Edif} and \eqref{dposX} one infers from
\eqref{Tjk1x} that
\begin{align}
%\label{sumTjk1}
\nonumber \biggl\|\sum_{k=1}^n T^{(1)}_{jk}x\biggr\| \leq &
[2\gamma_1+\gamma_2(d-c)]\,|\delta^{(m)}_j|
\left\|\bigl(\hsE(\lambda_j,d)-\hsE(\lambda_{j-1},d)\bigr)x\right\|.
\end{align}
Therefore
\begin{align}
\nonumber \biggl\|\sum_{j=1}^m\,\sum_{k=1}^n
T^{(1)}_{jk}x\biggr\|\leq& [2\gamma_1+\gamma_2(d-c)] \sum_{j=1}^m
|\delta^{(m)}_j|\,\left\|\bigl(\hsE(\lambda_j,d)-\hsE(\lambda_{j-1},d)\bigr)x\right\|\\
\nonumber \leq& [2\gamma_1+\gamma_2(d-c)]
\left(\sum_{j=1}^m|\delta^{(m)}_j|^2\right)^{1/2}
\left(\sum_{j=1}^m\left\|\bigl(\hsE(\lambda_j,d)-\hsE(\lambda_{j-1},d)
\bigr)x\right\|^2\right)^{1/2}\\
\nonumber \leq&[2\gamma_1+\gamma_2(d-c)] \,\left(\Max_{j=1}^m
|\delta^{(m)}_j|\right)^{1/2}
\left(\sum_{j=1}^m|\delta^{(m)}_j|\right)^{1/2}\\
\nonumber
&\qquad\times\left(\sum_{j=1}^m\bigl\langle\hsE(\lambda_j,d)-\hsE(\lambda_{j-1},d)x,
x\bigr\rangle\right)^{1/2}\\
\nonumber \leq&[2\gamma_1+\gamma_2(d-c)] \,\left(\Max_{j=1}^m
|\delta^{(m)}_j|\right)^{1/2}\sqrt{b-a}
\left(\bigl\langle\hsE(b,d)-\hsE(a,d)x,
x\bigr\rangle\right)^{1/2}\\
\label{Tjk1tot} \leq&[2\gamma_1+\gamma_2(d-c)] \,\left(\Max_{j=1}^m
|\delta^{(m)}_j|\right)^{1/2}\sqrt{b-a\,}\,\,
\|\hsE(b,d)-\hsE(a,d)\|\,\|x\|.
\end{align}

It is proven analogously that
\begin{align}
%\label{sumTjk2}
\nonumber
\biggl\|\sum_{j=1}^m T^{(2)}_{jk}x\biggr\| \leq &
[2\gamma_1+\gamma_2(b-a)]|\omega^{(n)}_k|
\left\|\bigl(\hsE(b,\mu_k)-\hsE(b,\mu_{k-1})\bigr)x\right\|
\end{align}
and then
\begin{align}
\label{Tjk2tot} \biggl\|\sum_{j=1}^m\,\sum_{k=1}^n
T^{(2)}_{jk}x\biggr\|\leq& [2\gamma_1+\gamma_2(b-a)]\left(\Max_{k=1}^n
|\omega^{(n)}_k|\right)^{1/2}\sqrt{d-c\,}\,\,
\|\hsE(b,d)-\hsE(b,c)\|\,\|x\|.
\end{align}

It only remains to find an estimate for the contribution to the
difference \eqref{JJdif} from the terms $T^{(3)}_{jk}$ given by
\eqref{T3jk}. By using identity \eqref{dpos1} and taking into
account \eqref{dLipschitz} it follows from \eqref{T3jk} that
\begin{align}
\nonumber \left\|T^{(3)}_{jk}x\right\|&\leq\gamma_2
\sum_{s=1}^{p_j-1}\,\sum_{t=1}^{q_k-1}
|\lambda_{j,s}-\lambda_{j,s-1}|\,|\mu_{k,t}-\mu_{k,t-1}|\,
\|\sE\bigl([\lambda_{j,s},\lambda_j)\times[\mu_{k,t},\mu_k)\bigr)x\|\\
\nonumber &\leq\gamma_2 \sum_{s=1}^{p_j-1}\,\sum_{t=1}^{q_k-1}
|\lambda_{j,s}-\lambda_{j,s-1}|\,|\mu_{k,t}-\mu_{k,t-1}|\,
\|\sE(\delta^{(m)}_j\times\omega^{(n)}_k)x\|\\
\nonumber
&\leq\gamma_2\,\,|\lambda_{j,p_j-1}-\lambda_{j-1}|\,|\mu_{k,q_k-1}-\mu_{k-1}|\,
\|\sE(\delta^{(m)}_j\times\omega^{(n)}_k)x\|\\
\label{T3jkx} &\leq \gamma_2\,\,|\delta^{(m)}_j|\,|\omega^{(n)}_k|
\|\sE(\delta^{(m)}_j\times\omega^{(n)}_k)x\|
=\gamma_2\,\,|\delta^{(m)}_j|\,|\omega^{(n)}_k|
\|\sE(\Delta^{(m,n)}_{jk})x\|\,.
\end{align}
Hence
\begin{align}
\nonumber \biggl\|\sum_{j=1}^m\sum_{k=1}^n
T^{(3)}_{jk}x\biggr\|&\leq
\sum_{j=1}^m\sum_{k=1}^n \left\|T^{(3)}_{jk}x\right\| \\
\nonumber &\leq \gamma_2\left(\sum_{j=1}^m\sum_{k=1}^n
|\delta^{(m)}_j|^2\,|\omega^{(n)}_k|^2\right)^{1/2}
\left(\sum_{j=1}^m\sum_{k=1}^n
\|\sE(\Delta^{(m,n)}_{jk})x\|^2\right)^{1/2} \\
\nonumber & \qquad =\gamma_2\left(\sum_{j=1}^m
|\delta^{(m)}_j|^2\right)^{1/2} \left(\sum_{k=1}^n
|\omega^{(n)}_k|^2\right)^{1/2} \left(\sum_{j=1}^m\sum_{k=1}^n
\langle \sE(\Delta^{(m,n)}_{jk})x,x\rangle\right)^{1/2}\\
\nonumber &\leq \gamma_2\left(\Max_{j=1}^m |\delta^{(m)}_j|\,\,
\Max_{k=1}^n |\omega^{(n)}_k|\right)^{1/2} \sqrt{(b-a)(d-c)}
\,\,\langle \sE(\Delta)x,x\rangle^{1/2} \\
\label{T3finest} &\leq \gamma_2\left(\Max_{j=1}^m |\delta^{(m)}_j|\,\,
\Max_{k=1}^n |\omega^{(n)}_k|\right)^{1/2} \sqrt{(b-a)(d-c)}
\|\sE(\Delta)\|\,\|x\|\,.
\end{align}
Now combining \eqref{JJdif}, \eqref{Tjk1tot}, \eqref{Tjk2tot}, and
\eqref{T3finest} one concludes that the integral sum \eqref{JmnL}
converges as $m\to\infty$, $n\to\infty$,
$\Max_{j=1}^m|\delta^{m}_j|\to0$, and
$\Max_{k=1}^n|\omega^{n}_k|\to0$  to a linear operator of
$\cB(\fK,\fH)$. The convergence takes place with respect to the
uniform operator topology.

The proof is complete.
\end{proof}
%%%%%%%%%%%%%%%%%%%%%%%%%%%%%%%%%%%%%%%%%%%%%%%%%%%%%%%%%%%%%%
\begin{remark}
It is an open problem whether or not the sufficient conditions
\eqref{Lipschitz}, \eqref{dLipschitz} for a function
$F:\Delta\to\cB(\fH,\fK)$ to be right-integrable on $\Delta$ with
respect to the spectral measure $dE(\lambda,\mu)$ in the sense of
the uniform operator topology are optimal. We also notice that
even in the case where the spectral family $\{E(\cdot)\}$ is
associated with a self-adjoint operator, a similar problem remains
unsolved. More precisely, it is known that if an operator-valued
function is Lipschitz on a finite interval in $\bbR$ then it is
integrable on this interval  in the sense of the uniform operator
topology  with respect to the spectral measure of any self-adjoint
operator (see \cite[Lemma 7.2 and Remark 7.3]{AdLMSr}). But, to
the best of our knowledge, it is still unknown whether the
requirement for the function to be Lipschitz is optimal.
\end{remark}
%%%%%%%%%%%%%%%%%%%%%%%%%%%%%%%%%%%%%%%%%%%%%%%%%%%%%%%%%%%%%%
\begin{corollary}
\label{Integr21} Assume that  for some $\gamma_1,\gamma_2>0$ an
operator-valued function $G:\,\Delta\rightarrow\cB(\fH,\fK)$
satisfies conditions \eqref{Lipschitz} and \eqref{dLipschitz}.
Then $G$ is left-integrable on $\Delta$ with respect to the
spectral measure $d\sE(\lambda,\mu)$ in the sense of the uniform
operator topology.
\end{corollary}
%%%%%%%%%%%%%%%%%%%%%%%%%%%%%%%%%%%%%%%%%%%%%%%%%%%%%%%%%%%%%%%%
\begin{proof}
The statement is an immediate consequence of Theorem \ref{Integr2} by
applying Lemma \ref{Integr1}.
\end{proof}

\begin{remark}
\label{remran} If the integral \eqref{defIntL} exists, its range
(as that of an operator of $\cB(\fH,\fK)$) lies in
$\fK_\Delta=\Ran\bigl(\sE(\Delta)\bigr)$.
\end{remark}

The next lemma gives a norm estimate for the integral
of an operator-valued function with respect to the spectral measure of
a normal operator.

\begin{lemma}
\label{Lnest}
Under the hypothesis of Theorem \ref{Integr2} the integral
\eqref{defIntR} satisfies the following norm estimate
\begin{align}
\nonumber
\Int_\Delta F(\lambda,\mu)\,d\sE(\lambda,\mu)\quad \leq\quad&
4\sup_{(\lambda,\mu)\in\Delta}\|F(\lambda,\mu)\| \\
\label{nestR}
& + 2\gamma_1\,(b-a+d-c)+\gamma_2(b-a)(d-c)].
\end{align}
\end{lemma}
\begin{proof}
Assume that $a=\lambda_0<\lambda_1<\ldots<\lambda_m=b$ and
$c=\mu_0<\mu_1<\ldots<\mu_n=b$. Let
$\delta_j^{(m)}=[\lambda_{j-1},\lambda_j)$, $j=1,2,\ldots,m$, and
$\omega_k^{(n)}=[\mu_{k-1},\mu_k)$, $k=1,2,\ldots,n$, be the
corresponding partitions of the intervals $[a,b)$ and $[c,d)$.
Also let $\Delta^{(m,n)}_{jk}=\delta_k^{(m)}\times\omega_k^{(n)}$.
Observe that the integral sum
\begin{equation}
\label{Jmn1} J_{m,n}=\sum_{j=1}^m \sum_{k=1}^n F(\lambda_{j-1},\mu_{k-1})
\sE(\Delta^{(m,n)}_{jk})
\end{equation}
reads
\begin{align*}
J_{m,n}=&\sum_{j=1}^m \sum_{k=1}^n
F(\lambda_{j-1},\mu_{k-1})
\bigl[\hsE(\lambda_j,\mu_k)-\hsE(\lambda_{j-1},\mu_k)
-\hsE(\lambda_{j},\mu_{k-1})+\hsE(\lambda_{j-1},\mu_{k-1})\bigr].
\end{align*}
By regrouping the terms and taking into account that $\lambda_{0}=a$,
$\lambda_{m}=b$, $\mu_0=c$, and $\mu_{n}=d$ one then verifies
that
\begin{align}
\label{JS}
J_{m,n}=S_1+S_2+S_3+S_4,
\end{align}
where
\begin{align*}
S_1=&F(\lambda_{m-1},\mu_{n-1})\hsE(b,d)
-F(a,\mu_{n-1})\hsE(a,d)
-F(\lambda_{m-1},c)\hsE(b,c)
+F(a,c)\hsE(a,c), \\
S_2=&\sum_{j=1}^{m-1}\left\{
\bigl[F(\lambda_{j-1},\mu_{n-1})-F(\lambda_j,\mu_{n-1})\bigr]
\hsE(\lambda_j,d)-
\bigl[F(\lambda_{j-1},c)-F(\lambda_j,c)\bigr]
\hsE(\lambda_j,c)\right\}, \\
S_3=&\sum_{k=1}^{n-1}\left\{
\bigl[F(\lambda_{m-1},\mu_{k-1})-F(\lambda_{m-1},\mu_k)\bigr]
\hsE(b,\mu_k)-
\bigl[F(a,\mu_{k-1})-F(a,\mu_k)\bigr]
\hsE(a,\mu_k)\right\},
\end{align*}
and
\begin{align*}
S_4=\sum_{j=1}^{m-1}\sum_{k=1}^{n-1}
\bigl[F(\lambda_{j},\mu_{k})-F(\lambda_{j-1},\mu_k)-
F(\lambda_j,\mu_{k-1})+F(\lambda_{j-1},\mu_{k-1})\bigr]
\hsE(\lambda_j,\mu_k).
\end{align*}
We have
\begin{equation}
\label{Ein}
\|\hsE(\lambda,\mu)\|\leq 1\quad\text{for \,any}\quad\lambda,\mu\in\bbR
\end{equation}
and hence
\begin{align}
\label{JS1}
S_1&\leq 4\sup_{(\lambda,\mu)\in\Delta}\|F(\lambda,\mu)\|.
\end{align}
By hypothesis the operator-valued function $F$ satisfies
estimates \eqref{Lipschitz} and \eqref{dLipschitz}. Taking into
account \eqref{Ein}, this implies
\begin{align}
\label{JS2}
S_2\leq & 2\gamma_1 \sum_{j=1}^{m-1}|\lambda_{j-1}-\lambda_j|<2\gamma_1(b-a),\\
\label{JS3}
S_3\leq & 2\gamma_1 \sum_{k=1}^{n-1}|\mu_{k-1}-\mu_k|<2\gamma_1(d-c),
\end{align}
and
\begin{align}
\label{JS4}
S_4\leq & \gamma_2 \sum_{j=1}^{m-1}\sum_{k=1}^{n-1}
|\lambda_{j-1}-\lambda_j||\mu_{k-1}-\mu_k|<\gamma_2(b-a)(d-c).
\end{align}
Combining \eqref{JS} and \eqref{JS1}--\eqref{JS4} and
passing in \eqref{Jmn1} to the limit as
$\Max_{j=1}^m |\delta_j^{(m)}|+\Max_{k=1}^n
|\omega_k^{(n)}|\rightarrow0$ one arrives at inequality \eqref{nestR}
which completes the proof.
\end{proof}

\medskip

Now let $\Omega=\{z\in\bbC\,|\,a\leq\Real z< b, c\leq\Img z< d\}$ be a
rectangle in $\bbC$ with finite $a$, $b$, $c$, and $d$. Assume
that $F:\Omega\to\cB(\fK,\fH)$ and
$G:\Omega\to\cB(\fH,\fK)$ are such that the corresponding
functions $F(\lambda+\ri \mu)$ and $G(\lambda+\ri\mu)$ of real
variables $\lambda\in[a,b)$ and $\mu\in[c,d)$ are resp. right- and
left-integrable with the
spectral measure $d\sE(\lambda,\mu)$ over the rectangle
$\Delta=[a,b)\times[c,d)$ in $\bbR^2$.
In this case we set
\begin{align}
\label{convr}
\int_\Omega F(z)dE(z)&=\int_\Delta F(\lambda+\ri\mu)d\sE(\lambda,\mu),\\
\label{convl}
\int_\Omega dE(z)G(z)&=\int_\Delta d\sE(\lambda,\mu)G(\lambda+\ri\mu),
\end{align}
where $dE(z)$ stands for the spectral measure of the same normal
operator as the measure $d\sE(\lambda,\mu)$ but on the complex plane.

Let $F:\Omega\to\cB(\fK,\fH)$ and $G:\Omega\to\cB(\fH,\fK)$ are
such that the $F(\lambda+\ri\mu)$ and $G(\lambda+\ri\mu)$,
$(\lambda,\mu)\in\Delta$, satisfy on $\Delta$ conditions
\eqref{Lipschitz} and \eqref{dLipschitz}. Then the integrals
\eqref{convr} and \eqref{convl} exist by Theorem \ref{Integr2} and
Corollary \ref{Integr21}, respectively. Clearly, by
\eqref{defIntR} and \eqref{defIntL} only the values of $F$ and $G$
on the support of the measure $d\sE(\lambda,\mu)$ contribute to
these integrals. Assuming that $dE(z)$ (or, equivalently,
$d\sE(\lambda,\mu)$) is the spectral measure associated with a
normal operator $C$,  by
\begin{align}
\label{convrC}
\int\limits_{\spec(C)\cap\Omega} F(z)dE(z)&=\int_\Omega F(z)dE(z),\\
\label{convlC}
\int\limits_{\spec(C)\cap\Omega} dE(z)G(z)&=\int_\Omega dE(z)G(z)
\end{align}
we define the corresponding integrals of $F$ and $G$ over the
part of the spectrum of $C$ lying in $\Omega$.

In particular, if a function $F:\spec(C)\to\cB(\fK,\fH)$ (resp.
$G:\spec(C)\to\cB(\fH,\fK)$) defined (only) on the spectrum of $C$
admits an extension $\widetilde{F}$ (resp. $\widetilde{G}$) to the
whole rectangle $\Omega$ in such a way that conditions
\eqref{Lipschitz} and \eqref{dLipschitz} hold for
$\widetilde{F}(\lambda+\ri\mu)$ (resp. for
$\widetilde{G}(\lambda+\ri\mu)$), $(\lambda,\mu)\in\Delta$, we set
\begin{align}
\label{convrCt}
\int\limits_{\spec(C)\cap\Omega} F(z)dE(z)&=\int_\Omega \widetilde{F}(z)dE(z),\\
\label{convlCt}
\int\limits_{\spec(C)\cap\Omega} dE(z)G(z)&=\int_\Omega dE(z)\widetilde{G}(z).
\end{align}
Clearly, the results in \eqref{convrCt} and \eqref{convlCt} do not
depend on the choice of the extensions $\widetilde{F}(z)$ and
$\widetilde{G}(z)$ (provided that $\widetilde{F}(\lambda+\ri\mu)$
and $\widetilde{G}(\lambda+\ri\mu)$ satisfy \eqref{Lipschitz} and
\eqref{dLipschitz}).

Finally, the improper weak, strong, or uniform integrals
\begin{equation}
\label{Imp} \Int_a^b \Int_c^d F(\lambda,\mu)\,d\sE(\lambda,\mu)
\text{\, and \,}
\Int_a^b \Int_c^d  d\sE(\lambda,\mu)\,G(\lambda,\mu)
\end{equation}
with infinite lower and/or upper bounds ($a=-\infty$ and/or
$b=+\infty$ and $c=-\infty$ and/or $d=+\infty$) are understood as
the limits, if they exist, of the integrals over
finite intervals in the corresponding topologies. For example,
$$
\Int_{-\infty}^\infty \Int_{-\infty}^\infty
d\sE(\lambda,\mu)\,G(\lambda,\mu)=
\mathop{\text{lim}}\limits_{\mbox{\scriptsize
$\begin{array}{c} a\downarrow -\infty, \,\,b\uparrow \infty \\
c\downarrow -\infty, \,\, d\uparrow \infty
\end{array}$}} \,\,
\Int_a^b \Int_c^d d\sE(\lambda,\mu)\,G(\lambda,\mu)\,.
$$

If $dE(z)$ (or,
equivalently, $d\sE(\lambda,\mu)$) is the spectral measure associated with a
normal operator $C$, we set
\begin{align*}
\Int_{\spec(C)} F(z)\,dE(z)&=
\Int_{-\infty}^{+\infty}\Int_{-\infty}^{+\infty}
\widetilde{F}(\lambda+\ri\mu)\,d\sE(\lambda,\mu),\\
\Int_{\spec(C)} dE(z)\,G(z)&=
\Int_{-\infty}^{+\infty}\Int_{-\infty}^{+\infty}
d\sE(\lambda,\mu)\,\widetilde{G}(\lambda+\ri\mu),
\end{align*}
assuming that $F:\spec(C)\to\cB(\fK,\fH)$ and
$G:\spec(C)\to\cB(\fH,\fK)$ admit extensions $\widetilde{F}$ and
$\widetilde{G}$ from $\spec(C)$ to the whole complex plane $\bbC^2$ in such
a way that
$\widetilde{F}(\lambda+\ri\mu)$ and
$\widetilde{G}(\lambda+\ri\mu)$ satisfy conditions
\eqref{Lipschitz}, \eqref{dLipschitz} as functions of the variables
$\lambda,\mu\in\bbR$.

We conclude this section with the the following natural
result.

\begin{lemma}
\label{CZero} Let $\fH$ and $\fK$ be Hilbert spaces. Assume that
$C$ is a (possibly unbounded) normal operator on $\fK$ and $dE(z)$
is the spectral measure associated with $C$. Let
$\Omega=\{z\in\bbC\,|\,a\leq\Real z< b, c\leq\Img z< d\}$ with
finite $a$, $b$, $c$, and $d$. Assume in addition that $G(z)$ is a
$\cB(\fH,\fK)$-valued function holomorphic on an open circle in
$\bbC$ containing the closure $\overline{\Omega}$ of the rectangle
$\Omega$. Then
\begin{equation}
\label{domincl}
\Ran\biggl(\int_\Omega dE(z)G(z)\biggr)\subset\Dom(C)
\end{equation}
and
\begin{equation}
\label{Zero}
\int_\Omega dE(z)\bigl(zG(z)\bigr)=C\int_\Omega dE(z)G(z).
\end{equation}
\end{lemma}
\begin{proof}
Inclusion \eqref{domincl} is proven by Remark \ref{remran} taking
into account definition \eqref{convl} of the right-side integral of
an operator-valued function with the measure $dE(z)$.

By hypothesis the functions $G(z)$ and $zG(z)$ are analytic on
$\Omega$. Hence, written in terms of the variables $\lambda=\Real
z\in[a,b)$, $\mu=\Img z\in[c,d)$, they
automatically satisfy conditions \eqref{Lipschitz} and
\eqref{dLipschitz}. Then by Theorem \eqref{Integr2} both
integrals
\begin{equation*}
\int_\Omega dE(z)G(z)\text{\, and \,}
\int_\Omega dE(z)\bigl(z\,G(z)\bigr)
\end{equation*}
exist in the sense of the uniform operator topology.

Assume that the open circle mentioned in hypothesis is centered at the
point $z_0$ and its radius equals $r$. Denote this circle by
$C_r(z_0)$ and write the Taylor formula
\begin{equation}
\label{Tay}
G(z)=\sum_{k=0}^{n} G_k (z-z_0)^k +R_n(z),\quad z\in C_r(z_0),
\end{equation}
where $G_k=\displaystyle\frac{G^{(k)}(z_0)}{k!}$ and $R_n(z)$ is
the remainder term. Since the function $G(z)$ is holomorphic on
$C_r(z_0)$ and the set $\Omega$ is a compact contained in
$C_r(z_0)$, the remainder $R_n(z)$ and its derivatives converge to
zero as $n\to\infty$ and the convergence is uniform with respect
to $z\in\Omega$.  In particular, this implies that
\begin{align}
\label{Se1}
&\qquad\qquad\sup\limits_{\lambda+\ri\mu\in\Omega}|S_n(\lambda,\mu)|
\mathop{\longrightarrow}\limits_{n\to\infty}0, \\
\label{Se2}
&\sup\limits_{\lambda+\ri\mu\in\Omega}\max\left\{\left|
\frac{\partial S_n(\lambda,\mu)}{\partial\lambda}\right|,\,
\left|\frac{\partial S_n(\lambda,\mu)}{\partial\mu}\right|
\right\}
\mathop{\longrightarrow}\limits_{n\to\infty}0,
\end{align}
and
\begin{equation}
\label{Se3}
\sup\limits_{\lambda+\ri\mu\in\Omega}\max\left\{\left|
\frac{\partial^2 S_n(\lambda,\mu)}{\partial\lambda^2}\right|,\,
\left|\frac{\partial^2 S_n(\lambda,\mu)}{\partial\mu^2}\right|,\,
\left|\frac{\partial^2
S_n(\lambda,\mu)}{\partial\lambda\partial\mu}\right|\right\}
\mathop{\longrightarrow}\limits_{n\to\infty}0,
\end{equation}
where
\begin{equation*}
S_n(\lambda,\mu)=zR_n(z)=(\lambda+\ri\mu)R_n(\lambda+\ri\mu),
\quad\lambda\in[a,b),\, \mu\in[c,d).
\end{equation*}
From \eqref{Se2} and \eqref{Se3} it follows that there is sequence
of $\gamma_n>0$, $n=0,1,2,\ldots,$ such that
\begin{equation}
\label{glim}
\gamma_n\to0\text{ \,  as \, }n\to\infty
\end{equation}
and
\begin{equation}
\label{LipS}
\begin{array}{c}
{\|S_n(\lambda,\mu)-S_n(\lambda',\mu')\|\leq
\gamma_n\,(|\lambda-\lambda'|+|\mu-\mu'|)},\\[0.5em]
{\|S_n(\lambda,\mu)-S_n(\lambda',\mu)-
S_n(\lambda,\mu')+S_n(\lambda',\mu')\|
\leq \gamma_n\,|\lambda-\lambda'||\mu-\mu'|}\\[0.5em]
\mbox{for any
$\lambda,\lambda'\in[a,b)$ and $\mu,\mu'\in[c,d)$.}
\end{array}
\end{equation}
By combining Corollary \ref{Integr21} and Lemma \ref{Lnest}, it
follows from \eqref{Se1}, \eqref{glim} and \eqref{LipS} that
\begin{equation}
\label{Dex1}
\biggl\|\int_\Omega dE(z)zR_n(z)\biggr\|
\mathop{\longrightarrow}\limits_{n\to\infty}0.
\end{equation}
By a similar reasoning one also infers that
\begin{equation*}
\bigl\|\int_\Omega dE(z)R_n(z)\bigr\|
\mathop{\longrightarrow}\limits_{n\to\infty}0
\end{equation*}
and hence
\begin{equation}
\label{Dex2}
\biggl\|C\int_\Omega dE(z)R_n(z)\biggr\|=
\biggl\|\bigl(CE(\Omega)\bigr)\int_\Omega dE(z)R_n(z)\biggr\|
\mathop{\longrightarrow}\limits_{n\to\infty}0
\end{equation}
since the product $CE(\Omega)$ is a bounded operator.

By \eqref{Tay} we have
\begin{align*}
  \int_\Omega dE(z)\bigl(z\,G(z)&=\sum_{k=0}^{n}
  \left(\int_\Omega  z(z-z_0)^k dE(z)\right) G_k +\int_\Omega dE(z)zR_n(z),
\end{align*}
which implies
\begin{align}
\nonumber
  \int_\Omega dE(z)\bigl(z\,G(z)&=\sum_{k=0}^{n}
  C(C-z_0)^k E(\Omega) G_k +\int_\Omega dE(z)zR_n(z)\\
\label{Cex1}
  &=C\sum_{k=0}^{n}
  \left(\int_\Omega  (z-z_0)^k dE(z)\right) G_k +\int_\Omega dE(z)zR_n(z),
\end{align}
by taking into account that $\Ran E(\Omega)\subset\Dom(C^l)$ for
any $l=1,2,\ldots$. On the other hand by \eqref{Tay} it follows that
\begin{align*}
%\label{Cex2}
C\int_\Omega dE(z)G(z)&=C\sum_{k=0}^{n}
  \left(\int_\Omega  (z-z_0)^k dE(z)\right) G_k +C\int_\Omega dE(z)R_n(z).
\end{align*}
Comparing this equality with \eqref{Cex1} yields
\begin{equation*}
 \int_\Omega dE(z)\bigl(z\,G(z)-C\int_\Omega dE(z)\,G(z)=
\int_\Omega dE(z)zR_n(z)-C\int_\Omega dE(z)R_n(z),
\end{equation*}
which by \eqref{Dex1} and \eqref{Dex2} completes the proof.
\end{proof}

\section{Norm of an operator with respect to the spectral measure}
\label{SEnorm}

The paper \cite{MM99} introduced the concept of the norm of a
bounded operator with respect the spectral measure of a
self-adjoint operator. This concept turned out to be a useful tool
in the study of the operator Sylvester and Riccati equations (see
\cite{AMM} for details). The goal of the present section is to
extend the above concept to the case where the spectral measure is
associated with a normal operator.

%%%%%%%%%%%%%%%%%%%%%%%%%%%%%%%%%%%%%%%%%%%%%%%%%%%%%%%%%%%%%%
\begin{definition}
\label{ENorm} Let $Y\in \cB(\fH,\fK)$ be a bounded operator from a
Hilbert space $\fH$ to a Hilbert space $\fK$ and let
$\{E(\Omega)\}_{\Omega\in\cA_\mathrm{Borel}(\bbC)}$ be the
spectral family of a (possibly unbounded) normal operator
on  $\fK$. Introduce
\begin{equation}
\label{enorma}
\|Y\|_{E}=\left(\Sup\limits_{\{\Omega_k\}}
\sum_k \|Y^*E(\Omega_k)Y\|\right)^{1/2},
\end{equation}
where the supremum is taken over finite (or countable) systems of
mutually disjoint Borel subsets $\Omega_k$ of the complex plane
$\bbC$, $\Omega_k\cap\Omega_l=\emptyset$, if $k\neq l$.  The
number $\|Y\|_{E}$ is called the norm of the operator $Y$ with
respect to the spectral measure $dE(z)$ or simply $E$-norm of $Y$.
For $Z\in\cB(\fK,\fH)$ the $E$-norm $\|Z\|_{E}$ is defined by
$\|Z\|_{E}=\|Z^*\|_{E}$.
\end{definition}

One easily verifies that if the norm $\|Y\|_{E}$ is finite then
$$
\|Y\|\le\|Y\|_{E}.
$$
If, in addition, $Y$ is a Hilbert-Schmidt operator, then
\begin{equation}
\label{hsch}
\|Y\|_{E}\le\|Y\|_2, \quad Y\in \cB_2(\fH, \fK),
\end{equation}
where $\|\cdot\|_2$ denotes the (Hilbert-Schmidt) norm on the ideal
$\cB_2(\fH,\fK)$ of Hilbert-Schmidt operators from $\fH$ to $\fK$.

The following statement is an extension of \cite[Lemma 10.7]{MM99}.

%%%%%%%%%%%%%%%%%%%%%%%%%%%%%%%%%%%%%%%%%%%%%%%%%%%%%%%%%%%%%%%%
\begin{lemma}
\label{IntegrN}
Assume that $C$ is a normal operator on the Hilbert space $\fK$.
Let an operator-valued function $F:\spec(C)\rightarrow\cB(\fH)$
be bounded
$$
\|F\|_\infty=\Sup\limits_{z\in\spec(C)}\|F(z)\|<\infty,
$$
and admit a bounded extension as a function of
$\lambda=\Real z$ and $\mu=\Img z$
from $\sigma(C)$ to the whole plane $\bbR^2$ which satisfies
conditions \eqref{Lipschitz} and \eqref{dLipschitz}. If the $E$-norm
$\|Y\|_{E}$ of an operator $Y\in \cB(\fH,\fK)$ with respect to
the spectral measure $dE(z)$ on $\bbC$ associated with $C$ is finite,
then the integrals
$$
\Int_{\spec(C)} dE(z)\,Y\,F(z) \quad \mbox{and} \quad
\Int_{\spec(C)} F(z)\,Y^*\,dE(z)
$$
exist in  the uniform operator  topology. Moreover, the following
bounds hold
\begin{equation}
\label{EstT1}
\biggl\|\Int_{\spec(C)} dE(z)\,Y\,F(z)\biggr\|
\leq\|Y\|_{E}\cdot\|F\|_\infty,
\end{equation}
\begin{equation}
\label{EstT2} \biggl\|\Int_{\spec(C)} F(z)\,Y^*\,dE(z)\biggr\|
\leq\|Y\|_{E}\cdot\|F\|_\infty.
\end{equation}
\end{lemma}
\begin{proof}
The proof is given for the case of the integral in \eqref{EstT1}.

For the extension of the function $F(\lambda+\ri\mu)$ from the set
$\sigma(C)$ to $\bbR^2$ we will use the notation $F(\lambda,\mu)$,
$\lambda,\mu\in\bbR$.

Let $[a,b)\subset\bbR$ and $[c,d)\subset\bbR$ be finite intervals,
$-\infty<a<b<\infty$, $-\infty<c<d<\infty$, and let
$\left\{\delta_j^{(m)}\right\}_{j=1}^{m}$ and
$\left\{\omega_k^{(n)}\right\}_{k=1}^{n}$ be partitions of $[a,b)$
and $[c,d)$, respectively. If
$\bigl(\delta_j^{(m)}\times\omega_k^{(n)}\bigr)\cap\sigma(C)\neq\emptyset$
choose $\xi_j\in\delta_j^{(m)}$ and $\zeta_k\in\omega_k^{(n)}$ in
such a way that $\xi_j+\ri\zeta_k\in\spec(C)$, that is,
$(\xi_j,\zeta_k)\in\sigma(C)$. If
$\bigl(\delta_j^{(m)}\times\omega_k^{(n)}\bigr)\cap\sigma(C)=\emptyset$
let $\xi_j$ and $\zeta_k$ be arbitrary points of the intervals
$\delta_j^{(m)}$ and $\omega_k^{(n)}$, respectively. Then taking
into account that
\begin{equation*}
\sE(\delta_j^{(m)}\times\omega_k^{(n)})\sE(\delta_s^{(m)}\times\omega_t^{(n)})=0
\text{\, if \,} j\neq s \text{\, or \,} k\neq t,
\end{equation*}
for any $x\in\fH$ one obtains
\begin{align}
\nonumber
&\left\|\sum_{j=1}^m\sum_{k=1}^n \sE(\delta_j^{(m)}\times\omega_k^{(n)})
Y F(\xi_j,\zeta_k)x\right\|^2\\
\nonumber
&\qquad=
\bigl\langle\sum_{j=1}^m\sum_{k=1}^n \sE(\delta_j^{(m)}
\times\omega_k^{(n)})
Y F(\xi_j,\zeta_k)x,
\sum_{s=1}^m\sum_{t=1}^n \sE(\delta_s^{(m)}\times\omega_t^{(n)})
Y F(\xi_s,\zeta_t)x\bigr\rangle\\
\nonumber
&\qquad=
\bigl\langle\sum_{j=1}^m\sum_{k=1}^n [F(\xi_j,\zeta_k)]^*Y^*\sE(\delta_j^{(m)}
\times\omega_k^{(n)})
Y F(\xi_j,\zeta_k)x,x\bigr\rangle\\
\nonumber
&\qquad=
\bigl\langle\sum_{j,k\,:\,\,
(\xi_j,\zeta_k)\in\sigma(C)} [F(\xi_j,\zeta_k)]^*Y^*\sE(\delta_j^{(m)}
\times\omega_k^{(n)})
Y F(\xi_j,\zeta_k)x,x\bigr\rangle.
\end{align}
Hence
\begin{align}
\nonumber
&\left\|\sum_{j=1}^m\sum_{k=1}^n \sE(\delta_j^{(m)}\times\omega_k^{(n)})
Y F(\xi_j,\zeta_k)x\right\|^2\\
\nonumber
&\qquad\qquad\leq \sup\limits_{(\xi,\zeta)\in\sigma(C)}\|F(\xi,\zeta)\|^2
\sum_{j,k\,:\,\,
(\xi_j,\zeta_k)\in\sigma(C)\cap\Delta} \|Y^*\sE(\delta_j^{(m)}\times\omega_k^{(n)})Y\|
\|x\|^2\\
\label{IntEst}
&\qquad\qquad\leq \|F\|_\infty^2 \|Y\|_{E,\Delta}^2
\|x\|^2,
\end{align}
where $\Delta=[a,b)\times[c,d)$ and
\begin{equation}
\label{enormar}
\|Y\|_{E,\Delta}=\left(\Sup\limits_{\{\Delta_k\}}
\sum_k \|Y^*\sE(\Delta_k)Y\|\right)^{1/2}.
\end{equation}
In \eqref{enormar} the supremum is taken over finite (or
countable) systems of Borel subsets $\Delta_k$ of the rectangle $\Delta$
such that $\Delta_k\cap\Delta_l=\emptyset$, if $k\neq l$. Obviously,
\begin{equation}
\label{enar0}
\|Y\|_{E,\Delta}\leq\|Y\|_E.
\end{equation}
and
\begin{equation}
\label{enar1}
\lim \limits_{
\begin{array}{c}
a\downarrow-\infty, \, b\uparrow\infty\\
c\downarrow-\infty, \, d\uparrow\infty
\end{array}}
\|Y\|_{E,\Delta}=\|Y\|_E.
\end{equation}

By hypothesis the function $F(\lambda,\mu)$ satisfies the assumptions of
Theorem \ref{Integr2} and hence by Corollary \ref{Integr21}
it is right-integrable on the rectangle $\Delta$ with respect
to the measure $d\sE(\lambda,\mu)$. From \eqref{IntEst} it follows that
\begin{align}
\label{Ifinest}
&\biggl\|\Int_{\Delta} d\sE(\lambda,\mu)\,Y\,F(\lambda,\mu)\biggr\|\leq
\|F\|_\infty\|Y\|_{E,\Delta}.
\end{align}
Thus \eqref{enar1} implies that
$$
\begin{array}{c}
\biggl\|\Int_{[a',b')\times[c',d')} d\sE(\lambda,\mu)\,Y\,F(\lambda,\mu)\biggr\|\rightarrow 0\\
\text{ as } a'\to\infty \text{ or } b'\to-\infty \,  \text{ \, (for } \, b'>a')  \text{ and/or }
 c'\to\infty \text{ or } d'\to-\infty \text{ \, (for } \, d'>c'),
\end{array}
$$
which together with \eqref{Ifinest} proves the existence of
the limit
$$
\Int_{\bbR^2} d\sE(\lambda,\mu)\,Y\,F(\lambda,\mu)=
\nlim \limits_{
\begin{array}{c}
a\downarrow-\infty, \, b\uparrow\infty\\
c\downarrow-\infty, \, d\uparrow\infty
\end{array}}
\Int_{[a,b)\times[c,d)} d\sE(\lambda,\mu)\,Y\,F(\lambda,\mu)
$$
and hence the convergence of the integral $\Int_{\spec(C)}
dE(z)\,Y\,F(z)$ in the sense of the operator norm topology  as
well as the bound \eqref{EstT1}. Then the convergence
of the integral $ \Int_{\spec(C)} F(z)\,Y^*\,dE(z)$ with respect to
the operator norm topology and the norm estimate \eqref{EstT2} are
proven by applying Lemma \ref{Integr1}.

The proof is complete.
\end{proof}

\section{Sylvester Equation}
\label{SSyl}

Assuming that the entry $C$ in the operator Sylvester equation
\eqref{SylInt} is a normal operator, the principal goal of this
section is to introduce a  Stieltjes integral representation for
the solution $X$ in terms of the spectral measure associated with
$C$.

We begin with recalling the concepts of weak, strong, and operator
solutions to the operator Sylvester equations.

\begin{definition}
\label{DefSolSyl}
Let $A$ and $C$ be densely defined possibly unbounded closed operators
on Hilbert spaces $\fH$ and $\fK$, respectively.  A bounded
operator $X\in \cB(\fH,\fK)$ is said to be a weak solution of
the Sylvester equation
\begin{equation}
\label{syl}
XA-CX=D, \quad D\in \cB(\fH,\fK),
\end{equation}
if
\begin{equation}
\label{sylw}
\lal XAf,g\ral-\lal Xf,C^* g\ral=\lal Df,g\ral
\quad \text{ for all } f\in \dom(A)
\text{ and } g\in \dom(C^*).
\end{equation}
A bounded operator $X\in \cB(\fH,\fK)$ is said to be a strong
solution of the Sylvester equation \eqref{syl} if
\begin{equation}
\label{ransyl}
\ran\bigl({X}|_{\dom(A)}\bigr)\subset\dom(C),
\end{equation}
and
\begin{equation}
\label{syls}
XAf-CXf=Df \quad \text{ for all } f\in \dom(A).
\end{equation}
Finally, a bounded operator $X\in \cB(\fH,\fK)$ is said to be
an operator solution of the Sylvester equation \eqref{syl} if
$$
\ran(X)\subset \dom (C),
$$
the operator $XA$ is bounded on $\dom (XA)=\dom(A)$, and the
equality
\begin{equation}
\label{sylext}
\overline{XA}-CX=D
\end{equation}
holds as an operator equality, where $\overline{XA}$
denotes the closure of
$XA$.
\end{definition}
%%%%%%%%%%%%%%%%%%%%%%%%%%%%%%%%%%%%%%%%%%%%%%%%%%%%%%%%%%%%%%%%%%%
Along with the Sylvester equation \eqref{syl} we also introduce
the dual equation
\begin{equation}
\label{sylZ}
YC^*-A^*Y=D^*,
\end{equation}
for which the notion of weak, strong, and operator solutions is
defined in a way analogous to that in Definition
\Ref{DefSolSyl}.
%%%%%%%%%%%%%%%%%%%%%%%%%%%%%%%%%%%%%%%%%%%%%%%%%%%%%%%%%%%%%%%%%%%

Clearly, if $X\in\cB(\fK,\fH)$ is an operator solution of the
Sylvester equation  \eqref{syl}, it is also a strong solution to
\eqref{syl}. In its turn, any strong solution is also a weak
solution. In fact, one does not need to distinguish between weak
and strong solutions to the Sylvester equation \eqref{syl} since
any weak solution to this equation is in fact a strong solution.

\begin{lemma}
\label{weak-strong}
Let $A$ and $C$ be densely defined possibly unbounded closed
operators on the Hilbert spaces $\fH$ and $\fK$, respectively.
If $X\in\cB(\fK,\fH)$ is a weak solution of the Sylvester equation
\eqref{syl} then $X$ is a strong solution of \eqref{syl}, too.
\end{lemma}
\begin{proof}
Given $f\in\Dom(A)$, introduce a linear functional $\mathfrak{l}_f$
with $\Dom(\mathfrak{l}_f)=\Dom(C^*)$ by
\begin{equation*}
\mathfrak{l}_f(g)=\lal C^*g,Xf\ral=\lal g,XAf\ral-
\lal g,Df\ral,\quad g\in\Dom(C^*).
\end{equation*}
Clearly, for any $f\in\Dom(A)$ the functional $\mathfrak{l}_f$
is bounded,
\begin{equation*}
|\mathfrak{l}_f(g)|=|\lal C^*g,Xf\ral|\leq c_f\|g\|,\quad g\in\Dom(C^*),
\end{equation*}
where $c_f=\|XAf\|+\|Df\|$. The functional $\mathfrak{l}_f$ is
also densely defined since $\Dom(\mathfrak{l}_f)=\Dom(C^*)$  is
dense in $\fK$ as domain of adjoint of a closed densely defined
operator. Thus, $Xf\in\Dom\bigl((C^*)^*\bigr)=\Dom(C)$ which
implies that \eqref{ransyl} holds and
\begin{equation*}
\lal g,XAf-CXf-Df\ral=0 \text{ \, for all \, }g\in\fH.
\end{equation*}
Hence \eqref{syls} also holds, completing the proof.
\end{proof}

It is easy to see that if one of the equations \eqref{syl} or
\eqref{sylZ} has a weak (and hence strong) solution then so does the other one.
%%%%%%%%%%%%%%%%%%%%%%%%%%%%%%%%%%%%%%%%%%%%%%%%%%%%%%%%%%%%%%%%%%%
\begin{lemma}
\label{XZsylw}
An operator $X\in\cB(\fH,\fK)$ is a weak (and hence strong)
solution to the Sylvester equation \eqref{syl} if and only if the
operator $Y=-X^*$ is a weak (and hence strong) solution to the
dual Sylvester equation \eqref{sylZ}.
\end{lemma}
\begin{proof}
The assertion is proven by combining Lemma \ref{weak-strong}
with \cite[Lemma 2.4]{AMM}.
\end{proof}

It is well known that if the spectra of the operators $A$ and $C$
are disjoint and one of them is a bounded operator then the
Sylvester equation $XA-CX=D$ has a unique solution.
This was first proven by by M.\,Krein in 1948.
Later, the same result was independently obtained by
Y.\,Daleckii \cite{D53} and M.\,Rosenblum \cite{R56}.
The precise statement is as follows.

\begin{lemma}
\label{Krein}
Let $A$ be a possibly unbounded densely defined closed operator
in the Hilbert space $\fH$ and $C$ a  bounded  operator on the
Hilbert space $\fK$ such that
$$
\spec(A)\cap\spec(C)=\emptyset
$$
and $D\in  \cB(\fH,\fK)$.  Then  the Sylvester equation
\eqref{syl} has a unique operator solution
\begin{equation}
\label{ESylKrein}
X=\frac{1}{2\pi\ri}
\int_{\Gamma} d\zeta\,(C-\zeta)^{-1}D (A-\zeta)^{-1},
\end{equation}
where $\Gamma$ is a union of closed contours in the complex plane
with total winding number $0$ around every point of $\spec(A)$ and
total winding number $1$ around every point of $\spec(C)$.
Moreover,
$$
\|X\|\le (2\pi)^{-1}|\Gamma|
\sup_{\zeta\in\Gamma} \|(C-\zeta)^{-1}\|\, \|(A-\zeta)^{-1}\|\,
\|D\|,
$$
where $|\Gamma|$ denotes the length of the contour $\Gamma$.
\end{lemma}
%%%%%%%%%%%%%%%%%%%%%%%%%%%%%%%%%%%%%%%%%%%%%%%%%%%%%%%%%%%%%%%%%%%

A relatively recent review of results on the Sylvester operator
equation \eqref{syl} with both bounded entries $A$ and $C$ and
applications of them to various problems can be found in
\cite{BR}.

If  $A$ and $C$ are unbounded densely defined closed operators,
even with separated spectra, then the Sylvester equation
\eqref{syl} may not have bounded solutions at all (a corresponding
example can be found in  \cite{P91}). Nevertheless, under some
additional assumptions equation \eqref{syl} is still solvable. A
review of known sufficient conditions for solvability of
\eqref{syl} in the case where both $A$ and $C$ are allowed to be
unbounded operators can be found in \cite[Section 2]{AMM}.

Now we prove the main result of this section:  if either $A$
or $C$ is normal, then a strong solution to the Sylvester
equation, if it exists, can be represented  in the form of an
operator Stieltjes integral. The corresponding representation
replaces the formula \eqref{ESylKrein} where the contour
integration might be impossible to perform in case of an unbounded
operator $C$.

%%%%%%%%%%%%%%%%%%%%%%%%%%%%%%%%%%%%%%%%%%%%%%%%%%%%%%%%%%%%%%%%%
\begin{theorem}
\label{SylUnb} Let $A$ be a possibly unbounded densely defined
closed operator on the Hilbert space $\fH$ and $C$ a normal
operator on the Hilbert space $\fK$ with the spectral family
$\{E_C(\Omega)\}_{\Omega\in\mathcal{A}_\mathrm{Borel}(\bbC)}$. Let
$D\in\cB(\fH,\fK)$ and suppose that $A$ and $C$ have disjoint
spectra, i.e.,
\begin{equation}\label{ACC}
\dist\bigl(\spec(A),\spec(C)\bigr)>0.
\end{equation}
Then the following statements are valid:
\begin{enumerate}

\item[{\rm$(i)$}] Assume that the Sylvester equation \eqref{syl}
has a strong solution $X\in\cB(\fH,\fK)$.  Then  $X$ is a unique
strong solution to \eqref{syl} and it can be represented in the
form of the Stieltjes integral
\begin{equation}
\label{EqX}
X=\int_{\spec(C)} dE_{C}(\zeta)
D(A-\zeta)^{-1},
\end{equation}
which converges in the sense of the strong operator
topology in $\cB(\fH,\fK)$.

Conversely, if the Stieltjes integral \eqref{EqX} exists in
the sense of the strong operator topology, then $X$ given by \eqref{EqX} is a
strong solution to \eqref{syl}.

\item[{\rm$(ii)$}] Assume that the dual Sylvester equation
\begin{equation}
\label{sylZs}
YC^*-A^*Y=D^*
\end{equation}
has a strong solution $Y\in\cB(\fK,\fH)$.  Then  $Y$ is a unique
strong solution to~\eqref{sylZs} and it can be
represented in the form of the Stieltjes operator integral
\begin{equation}
\label{EqXZ}
Y=-\int_{\spec(C^*)} (A^*-\zeta)^{-1} D^* E_{C^*}(d\zeta)
\end{equation}
that exists in the sense of the strong operator
topology in $\cB(\fK,\fH)$.
\end{enumerate}
\end{theorem}
%%%%%%%%%%%%%%%%%%%%%%%%%%%%%%%%%%%%%%%%%%%%%%%%%%%%%%%%%%%
\begin{proof} $(i)$ Set $\Delta=[a, b)\times[c,d)$
where $[a,b)$ and $[c,d)$ are finite real intervals.
Let $\{\delta_j\}$ be a finite system of mutually disjoint
intervals such that $[a,b)=\cup_j \delta_j$ and $\{\omega_k\}$
another finite system of mutually disjoint intervals such that
$[c,d)=\cup_k \omega_k$. Further, introduce the partition rectangles
$\Delta_{jk}=\delta_j\times\omega_k$. For the pairs $j,k$ such
that $\Delta_{jk}\cap\sigma(C)\ne \emptyset$ pick
$\zeta_{\Delta_{jk}}\in\bbC$ such that the point
$(\Real\zeta_{\Delta_{jk}},\Img\zeta_{\Delta_{jk}})\in\bbR^2$ belongs to
the intersection $\Delta_{jk}\cap\sigma(C)$. Applying to both
sides of \eqref{syls} the spectral projection
$\sE_{C}(\Delta_{jk})$, a short computation yields
\begin{equation}
\label{koki}
 \sE_{C}(\Delta_{jk})XAf-\zeta_{\Delta_{jk}} \sE_{C}(\Delta_{jk})Xf=
 \sE_{C}(\Delta_{jk})Df+\sE_{C}(\Delta_{jk})(C-\zeta_{\Delta_{jk}})Xf
\end{equation}
for any  $ f\in \dom(A)$.  Since
$(\Real\zeta_{\Delta_{jk}},\Img\zeta_{\Delta_{jk}})\in\Delta_{jk}\cap\sigma(C)$,
by \eqref{ACC} one concludes that $\zeta_{\Delta_{jk}}$ belongs to
the resolvent set of the operator $A$. Hence, \eqref{koki} implies
\begin{equation}\label{vagno}
 \sE_{C}(\Delta_{jk})X=\sE_{C}(\Delta_{jk})D(A-\zeta_{\Delta_{jk}})^{-1}
+(C-\zeta_{\Delta_{jk}})\sE_{C}(\Delta_{jk})X(A-\zeta_{\Delta_{jk}})^{-1}.
\end{equation}

Using
\eqref{vagno} one obtains
\begin{align}
\sum\limits_{j,k\,:\,\Delta_{j,k}\cap\,\sigma(C)\neq\emptyset}
\sE_{C}(\Delta_{j,k})X=&
\sum\limits_{j,k\,:\,\Delta_{j,k}\cap\,\sigma(C)\neq\emptyset}
\sE_{C}(\Delta_{j,k})D
(A-\zeta_{\Delta_{j,k}})^{-1}
\no \\
&+\sum\limits_{j,k\,:\,\Delta_{jk}\cap\,\sigma(C)\neq\emptyset}
(C-\zeta_{\Delta_{j,k}})\,\sE_{C}(\Delta_{j,k})
X(A-\zeta_{\Delta_{j,k}})^{-1}.
\label{dlinno}
\end{align}
The left hand side of \eqref{dlinno} can be computed explicitly:
\begin{equation}
\label{lili}
\sum\limits_{\Delta_{jk}\cap\,\sigma(C)\neq\emptyset}
\sE_{C}(\Delta_{jk})X
=\sE_{C}\bigl(\Delta\cap\,\sigma(C)\bigr)X=
\sE_{C}\bigl(\Delta\bigr)X=E_C(\Omega)X,
\end{equation}
where
\begin{equation}
\label{OmegaC}
\Omega=\{\zeta\in\bbC\,|\,\,a\leq\Real\zeta<b,\,\,c\leq\Img\zeta<d\}
\end{equation}
is the imbedding of the rectangle $\Delta$ into the complex plane
$\bbC$. Below we will also write the set $\Omega$ in the form
\begin{equation*}
\Omega=[a,b)\times\ri[c,d).
\end{equation*}

The first term on the right-hand side of \eqref{dlinno} is the
integral sum for the  Stieltjes integral  \eqref{EqX}.
More precisely, since $(A-\zeta)^{-1}$ is  analytic in a complex
neighborhood of $\Omega\cap\,\spec(C)$, by Theorem \Ref{Integr2}
and definitions \eqref{convlC}, \eqref{convlCt} one infers
\begin{align}
\label{Limit}
&{\nlim\limits_{\mathop{\rm max}\limits_{j,k}(|\delta_j|+|\omega_k|)\to 0}
\sum\limits_{j,k\,:\,\Delta_{j,k}\cap\,\sigma(C)\neq\emptyset}
\sE_{C}(\Delta_{j,k})D
(A-\zeta_{\Delta_{j,k}})^{-1}}\\
\nonumber
&\qquad=\int\limits_{\spec(C)\cap\Omega}
dE_{C}(\zeta)D(A-\zeta)^{-1}.
\end{align}
By the same reasoning for the last term on the right-hand side of
\eqref{dlinno} we have
\begin{align}
\nonumber
&\nlim\limits_{\mathop{\rm max}\limits_{j,k}(|\delta_j|+|\omega_k|)\to 0}
\sum\limits_{j,k\,:\,\Delta_{j,k}\cap\,\sigma(C)\neq\emptyset}
(C-\zeta_{\Delta_{j,k}})\,\sE_{C}(\Delta_{j,k})\,
X(A-\zeta_{\Delta_{j,k}})^{-1}\\
\label{putu}
&\qquad=C\int_{\spec(C)\cap\Omega}dE_C(\zeta)X(A-\zeta)^{-1}-
\int_{\spec(C)\cap\Omega}dE_C(\zeta)\,\zeta\,X(A-\zeta)^{-1}.
\end{align}
By using Lemma \ref{CZero} one easily proves that the right-hand side
of \eqref{putu} is zero, that is,
\begin{align}
\label{putu1}
&\nlim\limits_{\mathop{\rm max}\limits_{j,k}(|\delta_j|+|\omega_k|)\to 0}
\sum\limits_{j,k\,:\,\Delta_{j,k}\cap\,\sigma(C)\neq\emptyset}
(C-\zeta_{\Delta_{j,k}})\,\sE_{C}(\Delta_{j,k})\,
X(A-\zeta_{\Delta_{j,k}})^{-1}=0.
\end{align}

Passing to the limit  $\mathop{\rm
max}\limits_{j,k}(|\delta_j|+|\omega_k|)\to 0$ in \eqref{dlinno},
by combining \eqref{lili}, \eqref{Limit}, and \eqref{putu1} one concludes
that for any finite rectangle $\Omega$ of the form \eqref{OmegaC}
\begin{equation}
\label{semf}
E_{C}(\Omega)X=
\int\limits_{\spec(C)\cap\Omega}
dE_{C}(\zeta)D(A-\zeta)^{-1}.
\end{equation}
Since
\begin{equation*}
\slim\limits_{\mbox{\scriptsize$\begin{array}{c}
a,c\to-\infty\\b,d\to+\infty
\end{array}$}}
E_{C}\bigl([a,b)\times\ri[c,d)\bigr)X=X,
\end{equation*}
(\Ref{semf}) implies (\Ref{EqX}), which, in particular, proves
the  uniqueness of a strong solution to the Riccati equation
\eqref{syl}.

In order to prove the converse statement of $(i)$, assume that the
Stieltjes integral on the right-hand side of (\Ref{semf})
converges in the strong
operator topology  as $a,c\to-\infty$ and $b,d\to+\infty$ in
\eqref{OmegaC}. Denote the resulting integral by $X$.  Then,
(\Ref{semf}) holds for any finite rectangle $\Omega$ of the form
\eqref{OmegaC}.  This implies
that for any $f\in \dom(A)$ we have
\begin{align}
\nonumber
&{CE_{C}(\Omega)Xf-E_{C}(\Omega)XA f}\\
\nonumber
&\qquad=\int\limits_{\spec(C)\cap\Omega}
dE_{C}(\zeta)D(A-\zeta)^{-1}(\zeta-A)f\\
\nonumber
&\qquad=-\int\limits_{\spec(C)\cap\Omega}
dE_{C}(\zeta)Df=-E_{C}(\Omega)Df.
\end{align}
Hence,
\begin{align}
\nonumber
CE_{C}(\Omega)Xf&=E_{C}(\Omega)XA f-E_{C}(\Omega)Df\\
\label{predric}
&=E_{C}(\Omega)(XAf-Df) \\
\nonumber
& \quad \text{ \, for any \, } f\in
\dom(A).
\end{align}
In particular, \eqref{predric} implies that $ C E_{C}(\Omega)Xf$
with $\Omega$ given by \eqref{OmegaC} converges to $XAf-D f$ as
$a,c\to-\infty$ and $b,d\to+\infty$. Since
$\int_{\spec(C)\cap\Omega}\zeta dE_C(\zeta)Xf=CE_C(\Omega)Xf$,
from \eqref{predric} it also follows that
\begin{align*}
\int\limits_{\spec(C)\cap\Omega}|\zeta|^2 d\lal E_{C}(\zeta)Xf,Xf\ral
=\sup_{\Omega\in\cP} \bigl\|E_{C}(\Omega)(XA f-Df)\bigr\|^2
=\bigl\|XA f-Df\bigr\|^2
< \infty,
\end{align*}
where $\cP$ stands for the set of all rectangles in $\bbC$ of the form
\eqref{OmegaC}.
Hence,
\begin{equation}
\label{XfC}
Xf\in\dom(C).
\end{equation}
Then
\eqref{predric} can be rewritten as
\begin{equation}\label{BB}
 E_{C}(\Omega)C Xf
=E_{C}(\Omega)(XA f-D f) \text{ \, for any \, }\Omega\in\cP.
\end{equation}
Combining \eqref{XfC} and \eqref{BB} proves that $X$ is a strong
solution to the Sylvester equation \eqref{syl}.

%%%%%%%%%%%%%%%%%%%%%%%%%%%%%%%%%%%%%%%%%%%%%%%%%%%%%%%%%%%%%%%%%%%

{\rm$(ii)$}
Assume that the dual Sylvester equation~(\Ref{sylZ}) has a
strong solution $Y\in\cB(\fK,\fH)$. As  in the proof of $(i)$,
choose a finite rectangle $\Delta\subset\bbR^2$ such that
$\Delta \cap
\sigma({C^*})\ne \emptyset$.  Since
$\sE_{C^*}(\Delta)\fK\subset\dom({C^*})$, we have
$Y\sE_{C^*}(\Delta)f\in\dom(A^*)$ for any $f\in\fK$ by the definition
of a strong solution.  Take a point $\zeta_\Delta\in\spec({C^*})$ such
that $(\Real\zeta_\Delta,\Img\zeta_\Delta)\in \Delta$.
It follows from \eqref{ACC} that
$\zeta_\Delta\not\in\spec(A^*)$.  As in the proof of $(i)$, it is
easy to verify that
\begin{equation}
\label{vagnoZ}
\begin{array}{c}
Y\sE_{C^*}(\Delta)f=-(A^*-\zeta_\Delta)^{-1}D^*\sE_{C^*}(\Delta)f
-(A^*-\zeta_\Delta)^{-1}Y({C^*}-\zeta_\Delta)\sE_{C^*}(\Delta)f, \\
\end{array}
\end{equation}
which holds for any $f\in\fK.$

Next, let $[a,b)$ be a finite interval and $\{\delta_{j}\}$ a
finite system of mutually disjoint intervals such that
$[a,b)=\cup_j \delta_{j}$. Similarly, let $\{\omega_{k}\}$
be a finite system of mutually disjoint intervals partitioning
a finite interval $[c,d)$, i.e., $\cup_k \omega_{k}=[c,d)$.
Set $\Delta_{j,k}=\delta_j\times\omega_k$.
For $j,k$ such that $\Delta_{j,k}\cap
\sigma({C^*})\ne \emptyset$ pick a point
$\zeta_{\Delta_{j,k}}\in\spec({C^*})$ such that
$(\Real\zeta_{\Delta_{j,k}},\Img\zeta_{\Delta_{j,k}})\in\Delta_{j,k}$.
Using \eqref{vagnoZ} one then finds that
\begin{align}
\nonumber
Y\sE_{{C^*}}\bigl([a,b)\times[c,d)\bigr)f=&
-\sum\limits_{j,k:\,\Delta_{j,k}\cap\,\spec({C^*})\neq\emptyset}
(A^*-\zeta_{\Delta_{j,k}})^{-1}D^* \sE_{{C^*}}(\Delta_{j,k})f\\
\label{dlinnoZ}
&-\sum\limits_{j,k:\,\Delta_{j,k}\cap\,\spec({C^*})\neq\emptyset}
(A^*-\zeta_{\Delta_{j,k}})^{-1}Y\,({C^*}-\zeta_{\Delta_{j,k}})
\sE_{{C^*}}(\Delta_{j,k})f.
\end{align}
Equality \eqref{putu1} implies
\begin{equation}
\label{putuZ}
\nlim\limits_{\mathop{\rm max}\limits_{j,k}(|\delta_j|+|\omega_k|)\to 0}
\sum\limits_{j,k:\,\Delta_{j,k}\cap\,\spec({C^*})\neq\emptyset}
(A^*-\zeta_{\Delta_{j,k}})^{-1}Y({C^*}-\zeta_{\Delta_{j,k}})\sE_{{C^*}}(\Delta_{j,k})
=0.
\end{equation}
Thus, passing in (\Ref{dlinnoZ}) to the limit
as $\mathop{\rm max}\limits_{j,k}(|\delta_j|+|\omega_k|)\to 0$
one  infers that
\begin{equation}
\label{semfZ}
-\int\limits_{\spec({C^*})\cap\Omega}
(A^*-\zeta)^{-1}D^*dE_{{C^*}}(\zeta)f=Y\sE_{{C^*}}\bigl([a,b)\times[c,d)\bigr)f,
\end{equation}
where $\Omega$ is given by \eqref{OmegaC}. Since for any $f\in\fK$
\begin{equation*}
\mathop{\rm lim}\limits_{\mbox{\scriptsize$\begin{array}{c}
a,c\to-\infty\\b,d\to+\infty
\end{array}$}}
Y\sE_{{C^*}}\bigl([a,b)\times[c,d)\bigr)f=Y,
\end{equation*}
one concludes that the integral on the right-hand side of
\eqref{semfZ} converges as $a,c\to-\infty$ and $b,d\to+\infty$ in
the strong operator topology and thus (\Ref{EqXZ}) holds, which gives a
unique strong solution to the  dual Sylvester equation
\eqref{sylZs}.

%%%%%%%%%%%%%%%%%%%%%%%%%%%%%%%%%%%%%%%%%%%%%%%%%%%%%%%%%%%%%%%%%%%%%%%%

In order to prove   the converse statement of $(ii)$,  assume that
there exists the strong operator limit
\begin{equation}
\label{IntZab}
Y=\slim\limits_{\mbox{\scriptsize$\begin{array}{c}
a,c\to-\infty\\b,d\to+\infty
\end{array}$}}
\,\,\,\int\limits_{\spec({C^*})\cap([a,b)\times\ri[c,d))}
(A^*-\zeta)^{-1}D^*dE_{{C^*}}(\zeta) , \qquad Y\in\cB(\fK,\fH).
\end{equation}
Then for any finite $a$, $b$, $c$, and $d$ such that $a<b$ and
$c<d$ we have
\begin{equation}
\label{ZECab}
YE_{C^*}\bigl([a,b)\times\ri[c,d)\bigr)=
-\int\limits_{\spec({C^*})\cap([a,b)\times\ri[c,d))}
(A^*-\zeta)^{-1} D^* dE_{{C^*}}(\zeta).
\end{equation}
By \eqref{ACC} any point $\xi\in \spec({C^*})$ belongs to the
resolvent set of the operator $A$ and, hence, to the one of
$A^*$. Picking such a $\xi\in \spec({C^*})$, the
operator \eqref{ZECab} can be split into two parts
\begin{equation}
\label{ZEJJ}
YE_{{C^*}}\bigl([a,b)\times[c,d)\bigr)=J_1(a,b,c,d)+J_2(a,b,c,d),
\end{equation}
where
\begin{align}
\label{J1}
J_1(a,b,c,d)&=-(A^*-\xi)^{-1}D^* E_{{C^*}}([a,b)\times\ri[c,d)),\\
\label{J2}
J_2(a,b,c,d)&=+(A^*-\xi)^{-1}
\int\limits_{\spec({C^*})\cap ([a,b)\times\ri[c,d))}(\xi-\zeta)(A^*-\zeta)^{-1}D^*
dE_{{C^*}}(\zeta).
\end{align}
Using the functional calculus for the normal operator ${C^*}$
  one obtains
\begin{align*}
J_2(a,b,c,d)f=&
-(A^*-\xi)^{-1}\left(\int\limits_{\,\,\spec({C^*})\cap ([a,b)\times\ri[c,d))}
(A^*-\zeta)^{-1}D^* dE_{{C^*}}(\zeta)\right)({C^*}-\xi)f\\
&\text{ for any \,} f\in \dom ({C^*}).
\end{align*}
Thus, for $f\in \dom ({C^*})$ one concludes that
\begin{align*}
Yf=&
\lim\limits_{\mbox{\scriptsize$\begin{array}{c}
a,c\to-\infty\\b,d\to+\infty
\end{array}$}}YE_{{C^*}}\bigl([a,b)\times\ri[c,d)\bigr)f
\no \\
=&\lim\limits_{\mbox{\scriptsize$\begin{array}{c}
a,c\to-\infty\\b,d\to+\infty
\end{array}$}}
J_1(a,b,c,d)f
+
\lim\limits_{\mbox{\scriptsize$\begin{array}{c}
a,c\to-\infty\\b,d\to+\infty
\end{array}$}}J_2(a,b,c,d)f
\no \\
=&-(A^*-\xi)^{-1}D^*f
\no \\
&-(A^*-\xi)^{-1}\left(\int\limits_{\,\,\spec({C^*})}
(A^*-\zeta)^{-1}D^* dE_{{C^*}}(\zeta)\right)({C^*}-\xi)f
\no
\end{align*}
That is,
\begin{equation}
\label{ranZA}
Yf=-(A^*-\xi)^{-1}D^*f+(A^*-\xi)^{-1}Y({C^*}-\xi)f, \quad f\in \dom({C^*}),
\end{equation}
since
\begin{align}
\nonumber
&{\int_{\spec({C^*})}
(A^*-\zeta)^{-1}D^* dE_{{C^*}}(\zeta)} \\
\label{sIntDef}
&\qquad=\slim\limits_{\mbox{\scriptsize$\begin{array}{c}
a,c\to-\infty\\b,d\to+\infty
\end{array}$}}\,\,\,
\int\limits_{\spec({C^*})\cap ([a,b)\times\ri[c,d))}
(A^*-\zeta)^{-1}D^* dE_{{C^*}}(\zeta)=Y
\end{align}
by \eqref{IntZab}.
It follows  from~(\Ref{ranZA}) that
$Yf\in\dom(A^*)$ for any $f\in\dom({C^*})$ and, thus,
\begin{equation}
\label{ransylZ}
\ran\bigl({Y}|_{\dom({C^*})}\bigr)\subset\dom(A^*).
\end{equation}
Applying $A^*-\xi$ to both sides of the resulting equality
(\Ref{ranZA}) one infers that $Y$ is a strong solution to the
dual Sylvester equation~(\Ref{sylZs}), which completes the proof.
\end{proof}

%%%%%%%%%%%%%%%%%%%%%%%%%%%%%%%%%%%%%%%%%%%%%%%%%%%%%%%%%%
\begin{remark}
\label{CorXZ} Under the hypothesis of Theorem \ref{SylUnb} the
integral \eqref{EqX} converges in the sense of the strong operator
topology if and only if so does the integral \eqref{EqXZ}. This
can be seen by combining Theorem \ref{SylUnb} with Lemma
\ref{XZsylw}. The operators $X$ and $Y$ given by the integrals
\eqref{EqX} and \eqref{EqXZ} (if they exist in the sense of the
strong operator topology) are related to each other by $Y=-X^*$
\end{remark}
%%%%%%%%%%%%%%%%%%%%%%%%%%%%%%%%%%%%%%%%%%%%%%%%%%%%%%%%%%

\begin{lemma}
\label{CorXZ1} Assume that the hypothesis of Theorem \ref{SylUnb}
holds. If at least one of the integrals \eqref{EqX} and
\eqref{EqXZ} converges in the sense of the weak operator topology
then both of them converge  also in the sense of the strong
operator topology.
\end{lemma}

\begin{proof}
Suppose that the integral \eqref{EqX} converges in the sense of
the weak operator topology. Let
\begin{equation*}
X_\Omega=\int\limits_{\spec(C)\cap\Omega}
dE_{C}(\zeta)D(A-\zeta)^{-1},
\end{equation*}
where $\Omega=[a,b)\times\ri[c,d)$ is a finite rectangle in $\bbC$.
By the same reasoning as in the proof of equality \eqref{predric}
one obtains
\begin{equation}
\label{predsylw}
 X_\Omega Af-CX_\Omega f=E_{C}(\Omega)Df
\text{ \, for any \, }f\in\Dom(A),
\end{equation}
taking into account that $\Ran (X_\Omega)\subset\Ran
\bigl(E_C(\Omega)\bigr)$ by Remark \ref{remran} and hence
$X_\Omega f\in\Dom(C)$.
From \eqref{predsylw} it follows that
\begin{align}
\label{predsylw1}
&\lal X_\Omega Af,g\ral-\lal X_\Omega f,C^* g\ral=\lal E_C(\Omega)Df,g\ral\\
\nonumber
& \text{ \, for all } f\in \dom(A)
\text{ and } g\in \dom(C^*).
\end{align}
Passing in \eqref{predsylw1} to the limit as $a,c\to-\infty$ and
$b,d\to+\infty$ yields that $X$ given by the (improper) weak
integral \eqref{EqX} is a weak solution to the Sylvester equation
\eqref{syl} since by the assumption $\wlim X_{\Omega}=X$ and
since $\slim E_C(\Omega)=I$. Then by Lemma \ref{weak-strong} the
operator $X$ is a strong solution to \eqref{syl} and hence Theorem
\ref{SylUnb} $(ii)$ implies that the integral \eqref{EqX}
converges in the sense of the strong operator topology. By Remark
\ref{CorXZ} one concludes that so does the integral \eqref{EqXZ}.

Under the assumption that the integral
\eqref{EqXZ} converges in the sense of the weak operator, the
assertion is proven in a similar way.
\end{proof}

The next statement concerns sufficient conditions for the
existence of  a strong (and even operator) solution to the Sylvester
equation. The statement is an extension of \cite[Lemma 2.18]{AMM}.

%%%%%%%%%%%%%%%%%%%%%%%%%%%%%%%%%%%%%%%%%%%%%%%%%%%%%%%%%%
\begin{lemma}
\label{SylEC}
Assume hypothesis of Theorem \Ref{SylUnb}.
Suppose that the condition
\begin{equation}
\label{Res}
\Sup\limits_{\zeta\in\,\spec(C)}
\|(A-\zeta)^{-1}\|<\infty
\end{equation}
holds and the operator $D$ has a finite $E_C$-norm, that is,
\begin{equation}
\label{VarB}
\|D\|_{E_{C}}<\infty.
\end{equation}
Then equations \eqref{syl} and \eqref{sylZ} have unique strong
solutions $X\in\cB(\cH,\cK)$ given by \eqref{EqX} and
$Y\in\cB(\cK,\cH)$ given by \eqref{EqXZ}, respectively. Moreover,
$Y=-X^*$ and the Stieltjes integrals \eqref{EqX} and \eqref{EqXZ}
exist in the sense of the uniform operator topology.

Assume,  in addition, that
\begin{equation}
\label{ResHi}
\Sup\limits_{\zeta\in\,\spec(C)}
\|\zeta\,(A-\zeta)^{-1}\|<\infty.
\end{equation}
Then
\begin{align}
\label{XinC}
\ran(X)&\subset\dom(C),\\
\label{ZinA}
\ran(Y)&\subset\dom(A^*),
\end{align}
and thus $X$ and $Z$ appear to be operator solutions
to \eqref{syl} and \eqref{sylZ}, respectively.
\end{lemma}
%%%%%%%%%%%%%%%%%%%%%%%%%%%%%%%%%%%%%%%%%%%%%%%%%%%%%%%%%%%%%
We skip the proof since it almost literally repeats the proof of
Lemma 2.18 in \cite{AMM}. The only difference is in extending
the Stieltjes integration in the corresponding formulas from the
real axis to the complex plane.

\begin{remark}
Assume that
\begin{equation*}
\delta=\mathop{\rm dist}\bigl(\cW(A),\spec(C)\bigr)>0,
\end{equation*}
where $\cW(A)$ denotes the numerical range of $A$. Then by
\cite[Lemma V.6.1]{GK} it follows from \eqref{EqX} and
\eqref{VarB} that
\begin{equation}
\label{XYECr}
\|X\|_{E_C}\leq\frac{1}{\delta}\,\|D\|_{E_C},
\end{equation}
\end{remark}

%%%%%%%%%%%%%%%%%%%%%%%%%%%%%%%%%%%%%%%%%%%%%%%%%%%%%%%%%%%%%%%%%%%%
\begin{remark}
If $A$ is normal, then by \eqref{EqX} it immediately follows from
\eqref{VarB} that
\begin{equation}
\label{XYEC}
\|X\|_{E_C}\leq\frac{1}{d}\,\|D\|_{E_C},
\end{equation}
where $d=\mathop{\rm dist}\bigl(\spec(A),\spec(C)\bigr)$. In this
case one can also represent the operator $X$ in the form of a
double Stieltjes operator integral \cite{BS2003},
\begin{equation*}
X=\int_{\spec(C)}\int_{\spec(A)}
dE_{C}(\zeta)\,\dfrac{D}{z-\zeta}\,dE_{A}(z).
\end{equation*}
If $D\in\cB_2(\fH,\fK)$ then by \cite[Theorem 1]{BS67} the
operator $X$ is also Hilbert-Schmidt and the following estimate
holds
\begin{equation*}
\|X\|_2\leq\frac{1}{d}\,\|D\|_2,
\end{equation*}
which is sharp in the class of Hilbert-Schmidt operators.
\end{remark}

%%%%%%%%%%%%%%%%%%%%%%%%%%%%%%%%%%%%%%%%%%%%%%%%%%%%%%%%%%%%%%%%

%%%%%%%%%%%%%%%%%%%%%%%%%%%%%%%%%%%%%%%%%%%%%%%%%%%%%%%%%%%%%%%%
\section{Riccati equation}
\label{SecRic}
%%%%%%%%%%%%%%%%%%%%%%%%%%%%%%%%%%%%%%%%%%%%%%%%%%%%%%%%%%%%%%%%

There are at least three approaches that allow to tackle the
Riccati equations involving operators on infinite-dimensional
Hilbert spaces. The first of these approaches, going back to
C.\,Da\-vis \cite{Davis:58} and Halmos \cite{Halmos:69}, is based
on a deep connection between theory of Riccati equations and
results on variation of invariant subspaces of an operator under
perturbation. We refer to the recent publication \cite{KMM2}
discussing this purely geometric approach and its present status
in great detail. Here we only mention that such an approach is
essentially restricted to the operator Riccati equations
associated with self-adjoint block operator matrices, that is, to
the case of \eqref{RicInt} with $A=A^*$, $C=C^*$, and $D=B^*$.
Notice that the sharp norm estimates for variation of spectral
subspaces under a perturbation obtained in \cite{DK70,KMM4,MotSel}
imply the corresponding sharp norm estimates for solutions of
the associated Riccati equations.

The other approach is based on the factorization theorems for
holomorphic operator-valued functions proven by Markus and Matsaev
\cite{MarkusMatsaev75}  and by Virozub and Matsaev
\cite{VirozubMatsaev}. Several existence results for operator
Riccati equations have been obtained within this  approach (see
\cite{KMM3,MenShk}) including an existence result \cite{LMMT} for
the case where the entries $A$ and $C$ are allowed to be
non-self-adjoint operators.

The third approach \cite{AMM,MM99,MotSPbWorkshop} (see also
\cite{AdLT} and \cite{HMM}) is closely related to the integral
representation \eqref{EqX} for the solution of the operator
Sylvester equation in the form of an operator integral. Using this
representation allows one to rewrite the Riccati equation in the
form of an equivalent integral equation that admits an application
of Banach's Fixed Point Principle.  So far, only the Riccati
equations \eqref{RicInt} with at least one of the entries $A$ and
$C$ being a self-adjoint operator were studied within such an
approach. In this section we derive consequences of the integral
representation \eqref{EqX} that work for more general Riccati
equations \eqref{RicInt}, where one of the entries $A$ and $C$ is
merely a normal operator.

First, we recall the concepts of weak, strong, and operator
solutions to the operator Riccati equations.

%%%%%%%%%%%%%%%%%%%%%%%%%%%%%%%%%%%%%%%%%%%%%
\begin{definition}
\label{DefSolRic} Assume that  $A$ and $C$ are  possibly unbounded
densely defined closed operators on the Hilbert spaces $\fH$ and
$\fK$, respectively. Let $B$ and $D$ be bounded operators from $\fK$
to $\fH$ and from $\fH$ to $\fK$, respectively.

%Let $B\in\cB(\fK,\fH)$
%and $D\in\cB(\fH,\fK)$.

A bounded operator $X\in\cB(\fH,\fK)$ is
said to be a weak solution of the Riccati equation
\begin{equation}
\label{RicABCD}
XA-CX+XBX=D
\end{equation}
if
$$
%%%\label{Ricweek}
%
\begin{array}{c}
\lal XAf,g\ral-\lal Xf,C^*g\ral+\lal XBXf,g\ral=\lal Df,g\ral\\[0.5em]
 \text{ for all } f\in \dom (A)
\text{ and } g\in \dom(C^*).
\end{array}
$$

A bounded operator $X\in\cB(\fH,\fK)$ is called a strong
solution of the Riccati equation \eqref{RicABCD} if
\begin{equation}
\label{ranric}
\ran\bigl({X}|_{\dom(A)}\bigr)\subset\dom(C),
\end{equation}
and
\begin{equation}
\label{rics}
XAf-CXf+XBXf=Df  \text{ \, for all \, } f\in \dom(A).
\end{equation}

Finally,  $X\in \cB(\fH,\fK)$ is said to be
an operator solution of the Riccati equation \eqref{RicABCD} if
$$
\ran(X)\subset \dom (C),
$$
the operator $XA$ is bounded on $\dom(XA)=\dom(A)$,
and equality
\begin{equation}
\label{Ricext}
\overline{XA}-CX+XBX=D
\end{equation}
holds as an operator equality, where $\overline{XA}$ denotes the closure of
$XA$.
\end{definition}
%%%%%%%%%%%%%%%%%%%%%%%%%%%%%%%%%%%%%%%%%%%%%%%%%%%%%%%%%%%%%%%%%%%

Along with the Riccati equation \eqref{RicABCD} we also introduce
the dual equation
\begin{equation}
\label{RicK}
YC^* - A^*Y + YB^*Y=D^*.
\end{equation}

The following assertion is a corollary of Lemma \ref{weak-strong}.

\begin{lemma}
\label{weak-strongR} Let $A$ and $C$ be densely defined possibly
unbounded closed operators on the Hilbert spaces $\fH$ and $\fK$,
respectively. If $X\in\cB(\fK,\fH)$ is a weak solution of the
Riccati equation \eqref{RicABCD} then $X$ is also a strong
solution of \eqref{RicABCD}.
\end{lemma}
\begin{proof}
The assumption that $X$ is a weak solution to the Riccati equation
\eqref{RicABCD} implies that $X$ is a weak solution to the Sylvester
equation
\begin{equation}
\label{Sylaux}
X\widetilde{A}-CX=D,
\end{equation}
where
\begin{equation}
\label{hatA}
\widetilde{A}=A+BX \text{ \, with \, } \Dom(\widetilde{A})=\Dom(A)
\end{equation}
is a closed densely defined operator on $\fH$. Hence by Lemma
\ref{weak-strong} the operator $X$ is also a strong solution to
\eqref{Sylaux}, that is,
$\Ran(X|_{\Dom(\widetilde{A})})\subset\Dom(C)$ and
\begin{equation*}
X\widetilde{A}f-CXf=Df \text{ \, for all \, } f\in\Dom(\widetilde{A}).
\end{equation*}
Taking into account \eqref{hatA}, one then concludes  that $X$ is
a strong solution to \eqref{RicABCD}.
\end{proof}

The next statement is a direct corollary of Lemma \ref{XZsylw}.

%%%%%%%%%%%%%%%%%%%%%%%%%%%%%%%%%%%%%%%%%%%%%%%%%%%%%%%%%%%%%%%%
\begin{lemma}
\label{QKricw}
Let $A$ and $C$ be densely defined possibly unbounded
closed operators on the Hilbert spaces $\fH$ and $\fK$,
respectively, and $B\in\cB(\fK,\fH)$, $D\in\cB(\fH,\fK)$.
Then $X\in\cB(\fH,\fK)$ is a weak (and hence strong) solution to the Riccati equation
\eqref{RicABCD} if and only if $Y=-X^*$ is a weak (and hence strong)
solution to the dual Riccati equation \eqref{RicK}.
\end{lemma}
%%%%%%%%%%%%%%%%%%%%%%%%%%%%%%%%%%%%%%%%%%%%%%%%%%%%%%%%%%%%%%%%

Throughout the remaining part of the section we assume the
following hypothesis.
%%%%%%%%%%%%%%%%%%%%%%%%%%%%%%%%%%%%%%%%%%%%%%%%%%%%%%%%%%%%%%%%
\begin{hypothesis}
\label{hhqq}
Assume that $\fH$ and $\fK$ are Hilbert spaces, $A$ is a possibly
unbounded densely defined closed operator on $\fH$ and $C$
a normal operator on $\fK$.  Also assume that
$B\in\cB(\fK,\fH)$ and $D\in\cB(\fH,\fK)$.
\end{hypothesis}
%%%%%%%%%%%%%%%%%%%%%%%%%%%%%%%%%%%%%%%%%%%%%%%%%%%%%%%%%%%%%%%%

The representation theorems of  Sec. \ref{SSyl} for solutions
of the Sylvester equation provide us with iteration schemes
allowing one to prove solvability of the Riccati equations by
using fixed point theorems.

\begin{theorem}
\label{basics}
Assume Hypothesis \Ref{hhqq}.  Then the following statements
hold.

\begin{enumerate}

\item[{\rm$(i)$}] Assume, in addition to Hypothesis \Ref{hhqq}, that
\begin{equation}
\label{ABQCsep}
\dist\bigl(\spec(A+BX),\spec(C)\bigr)>0.
\end{equation}
Then $X\in\cB(\fH,\fK)$ is a weak (and hence strong) solution to
the Riccati equation \eqref{RicABCD} if and only if $X$ is a
solution of the equation
\begin{equation}
\label{EqQ}
X=\int_{\spec(C)} dE_{C}(\zeta)
D\,(A+BX-\zeta)^{-1},
\end{equation}
where  the operator Stieltjes integral exists in the sense of
the weak (and hence strong) operator topology in $\cB(\fH,\fK)$.

\item[{\rm$(ii)$}]
Assume, in addition to Hypothesis \Ref{hhqq}, that
$Y\in\cB(\fK,\fH)$ and
\begin{equation}
\label{ABKCsep}
\dist\bigl(\spec(A^*-YB^*),\spec(C^*)\bigr)>0.
\end{equation}
Then the operator $Y$ is a weak (and hence strong)  solution to
the dual Riccati equation \eqref{RicK} if and only if $Y$
satisfies the equation
\begin{equation}
\label{BasicK}
Y=-\int_{\spec(C^*)} (A^*-YB^*-\zeta)^{-1}D^*dE_{C^*}(\zeta),
\end{equation}
where the operator Stieltjes integral exists in the sense of the
weak (and hence strong) operator topology.

\end{enumerate}
\end{theorem}
%%%%%%%%%%%%%%%%%%%%%%%%%%%%%%%%%%%%%%%%%%%%%%%%%%%%%%%%%%%%%%%%%%%%%%
\begin{proof} $(i)$ The operator $X$ is a weak (and hence strong)
solution to \eqref{RicABCD} if and only if $X$ is a weak (and
hence strong) solution to the equation
$$
X\tA-CX=D,
$$
where
$$
\tA=A+BX.
$$
Applying Theorem \Ref{SylUnb} $(i)$ and Lemma \ref{CorXZ1}
completes the proof of $(i)$.

$(ii)$
 The operator $Y$ is a weak (and hence strong) solution to
\eqref{RicK} if and only if $Y$ is a weak (and hence strong) solution to the
equation
$$
YC-\widehat AY=D^*,
$$
where
$$
\widehat A=A-YB^*.
$$
Applying Theorem \Ref{SylUnb} $(ii)$ and Lemma \ref{CorXZ1}
completes the proof of $(ii)$.

The proof is complete.
\end{proof}
%%%%%%%%%%%%%%%%%%%%%%%%%%%%%%%%%%%%%%%%%%%%%%%%%%%%%%%%%%%%%%%%%%%%%%

The following statement is a direct corollary of Lemma
\Ref{SylEC}.
%%%%%%%%%%%%%%%%%%%%%%%%%%%%%%%%%%%%%%%%%%%%%%%%%%%%%%%%%%%%%%%%
\begin{theorem}
\label{RicEC}
Assume Hypothesis \Ref{hhqq} and
let  $D$ have  a finite norm with respect to
the spectral measure of the normal operator $C$, that is,
\begin{equation}
\label{VarDC}
\|D\|_{E_{C}}<\infty.
\end{equation}
Assume, in addition, that a bounded operator $X$ from $\fH$ to $\fK$
is a weak solution of the Riccati equation \eqref{RicABCD} such that
\begin{equation}
\label{ABQCsep1}
\dist\bigl(\spec(A+BX),\spec(C)\bigr)>0,
\end{equation}
and  that the condition
\begin{equation}
\label{Ress}
\Sup\limits_{\zeta\in\,\spec(C)}
\|(A+BX-\zeta)^{-1}\|<\infty
\end{equation}
holds.

Then  $X$ is a strong solution to
\eqref{RicABCD} and  the operator $Y=-X^*$
is a strong solution to the dual Riccati equation
\eqref{RicK}.

 The strong solutions $X$ and $Y$ admit the representations
 \begin{align}
\label{EqQ1}
X&=\int_{\spec(C)} dE_{C}(\zeta)
D\,(A+BX-\zeta)^{-1}, \\
\label{BasicK1}
Y&=-\int_{\spec(C)} (A-YB^*-\zeta)^{-1}D^*dE_{C^*}(\zeta),
\end{align}
where the operator Stieltjes integrals
 exist in the sense of the uniform operator
topology.
 Hence, the operators $X$ and $Y$ have finite
$E_C$--norm and the following bound holds true
\begin{equation}
\label{QEC0}
\|Y\|_{E_C}=\|X\|_{E_C}\leq \|D\|_{E_C}
\Sup_{\zeta\in\spec(C)}\|(A+BX-\zeta)^{-1}\|.
\end{equation}
If, in this case, instead of \eqref{Ress} the following condition
holds
\begin{equation}
\label{ResHiQ}
\Sup\limits_{\zeta\in\,\spec(C)}
\|\zeta\,(A+BX-\zeta)^{-1}\|<\infty,
\end{equation}
then
\begin{equation*}
\ran(X)\subset\dom(C), \quad \ran(Y)\subset\dom(A^*)
\end{equation*}
and, hence, the strong solutions $X$ and $Y$ appear to be operator
solutions to the Riccati equations \eqref{RicABCD} and
\eqref{RicK}, respectively.
\end{theorem}
%%%%%%%%%%%%%%%%%%%%%%%%%%%%%%%%%%%%%%%%%%%%%%%%%%%%%%%%%%%

In the case where the spectrum of the normal operator $C$ is
separated from the numerical range $\cW(A)$ of the operator $A$ we
are able to prove the existence of a fixed point for the mapping
\eqref{EqQ}, provided that the operators  $B$ and $D$ satisfy
certain ``smallness" assumptions. If, in addition, $A$ is also
a normal operator, we prove the existence of such a fixed point
under weaker assumptions.

%%%%%%%%%%%%%%%%%%%%%%%%%%%%%%%%%%%%%%%%%%%%%%%%%
\begin{theorem}
\label{QsolvN}
Under Hypothesis \Ref{hhqq} assume that $B\neq0$ and either
\begin{equation}
\label{tsemib}
(i)\,\, d=\dist\bigl(\cW(A),\spec(C)\bigr)>0
\qquad\quad\qquad\qquad\qquad\qquad\qquad\qquad\qquad\,\,
\end{equation}
or

$(ii)$ $A$ is a normal operator on $\fH$ and
\begin{equation}
\label{tsemibd}
  d=\dist\bigl(\spec(A),\spec(C)\bigr)>0.
\end{equation}
Also assume that the operator $D$ has a finite
$E_C$--norm and
\begin{equation}
\label{BDest}
\sqrt{\|B\|\,\|D\|_{E_C}}<\frac{d}{2}.
\end{equation}
Then the Riccati equation \eqref{RicABCD} has a unique strong solution
in the ball
\begin{equation}
\label{QEst2}
\left\{X\in\cB(\fH,\fK)\,\big|\,\,\,
\|X\|<\|B\|^{-1}\left(d-\sqrt{\|B\|\,\|D\|_{E_C}}\right)\right\}.
\end{equation}
Moreover, the strong solution $X$ has a finite $E_C$--norm that satisfies
the bound
\begin{equation}
\label{QEst1}
\|X\|_{E_C}\leq
 \frac{1}{\|B\|}\,
\left(\frac{d}{2}-\sqrt{\frac{d^2}{4}-\|B\|\,\|D\|_{E_C}}\right).
\end{equation}
In particular, for
\begin{equation}
\label{BDestSum}
\|B\|<\dfrac{d}{2} \text{ \, and \, } \|B\|+\|D\|_{E_C}<d
\end{equation}
or
\begin{equation}
\label{BDestSum1}
\|B\|\geq\dfrac{d}{2} \, \, \, \left(\text{and \, }
\|D\|_{E_C}<\dfrac{d^2}{4\|B\|}\right)
\end{equation}
% BYLO: \begin{equation}
%\label{BDestSum}
%\|B\|+\|D\|_{E_C}<d,
%\end{equation}
the strong solution  $X$ is a strict contraction in
both the $E_C$-norm and uniform operator topologies,
\begin{equation}
\label{Xstrict}
\|X\|\le \|X\|_{E_C}<1.
\end{equation}
\end{theorem}
%%%%%%%%%%%%%%%%%%%%%%%%%%%%%%%%%%%%%%%%%%%%%%%%%%
\begin{proof}
Technically, the proof of the assertion is very similar to that of
the second part of Theorem 3.6 in \cite{AMM} and it follows the same
line in both cases $(i)$ and $(ii)$. Therefore we present here only the
proof for the case~$(i)$.

Assume the hypothesis with assumption $(i)$.

We notice, first, that under this assumption the
resolvent $(A-\zeta)^{-1}$ for $\zeta\in\spec(C)$ is a bounded
operator since $\spec(A)\subset\overline{\cW(A)}$. Moreover, by
\cite[Lemma V.6.1]{GK} the following inequality holds
\begin{equation*}
\|(A-\zeta)^{-1}\|\leq\dfrac{1}{\dist(\zeta,\cW(A))},
\end{equation*}
which  \eqref{tsemib} yields
\begin{equation}
\label{Enew1}
\sup_{\zeta\in\spec(C)}\|(A-\zeta)^{-1}\|\leq\frac{1}{d}.
\end{equation}

Given $r\in (0,{d}/\|B\|)$, denote by $\cO_r$ the
closed $r$-neighbourhood of the zero operator in $\cB(\fH,\fK)$,
i.e. $\cO_r=\{X\in\cB(\fH,\fK)\,|\,\, \|X\|\leq r\}$.
Clearly, for any $X\in\cO_r$ by \eqref{Enew1} we have
\begin{align}
\nonumber
\sup_{\zeta\in\spec(C)}\left\|\left(I+(A-\zeta)^{-1}BX\right)^{-1}\right\|
&\leq\sup_{\zeta\in \spec(C)}\frac{1}{1-\|(A-\zeta)^{-1}\|\,\|B\|\,\|X\|}\\
\label{Enew2}
&\leq\dfrac{1}{1-\dfrac{\|B\| r}{{d}}}.
\end{align}
Therefore, the inverse operators involved in the identity
\begin{align}
\label{ResIdentity}
(A+BX-\zeta)^{-1}=&
\left(I+(A-\zeta)^{-1}BX\right)^{-1}(A-\zeta)^{-1},\\
\nonumber
&\zeta\in \spec(C),
\,\,X\in \cO_r,
\end{align}
are well defined. By \eqref{Enew1} and \eqref{Enew2} from
\eqref{ResIdentity} it follows that
\begin{align}
\label{ResExt}
\sup_{\zeta\in \spec(C)}\|(A+BX-\zeta)^{-1}\| \leq&
\displaystyle\frac{1}{{d}-\|B\| r},
\text{ whenever
$ X\in \cO_r.$}
\end{align}

Since \eqref{ResExt} holds and the operator $D$ has a finite
$E_C$-norm,
the mapping
\begin{equation}
\label{FX}
F(X)=\int_{\spec(C)} dE_{C}(\zeta)D\,(A+BX-\zeta)^{-1},
\end{equation}
from $\cB(\fH)$ to $\cB(\fK)$ is well defined on $\dom(F)=\cO_r.$
Notice that by Lemma \ref{IntegrN} the operator Stieltjes integral
in \eqref{FX} exists even in the sense of the uniform operator
topolgy.

The assumption that $X\in\cO_r$ implies that the numerical range
$\cW(A+BX)$ of the operator $A+BX$ lies in the closed
$(\|B\|r)$-neighborhood of $\cW(A)$. Since $\|B\|r<d$ and
$\spec(A+BX)\subset\overline{\cW(A+BX)}$, one then concludes that
$\dist\bigl(\spec(A+BX),\spec(C)\bigr)\geq d-\|B\|r>0$. Hence from
Theorem \Ref{basics} $(i)$ it follows that any fixed point of the
mapping $F$ in the ball $\cO_r$, if it exists, appears to be a
strong solution to the Riccati equation \eqref{RicABCD}.

Using \eqref{ResExt}
we obtain the following two estimates
\begin{align}
\|F(X)\|\leq\|F(X)\|_{E_C}&\leq \|D\|_{E_C}
\Sup_{\zeta\in\spec(C)}\|(A+BX-\zeta)^{-1}\|
\no\\
&\le\frac{\|D\|_{E_C}}{{d}-\|B\|r}, \quad X\in \cO_r,
\label{QQFC}
\end{align}
and
\begin{align}
\no
&\|F(X_1)-F(X_2)\|\\
&\qquad\leq\|F(X_1)-F(X_2)\|_{E_C}
\no \\
&\quad\qquad=\left\|\,\int_{\spec(C)} dE_{C}(\zeta)D\,
 (A+BX_1-\zeta)^{-1}\,B(X_2-X_1)\,(A+BX_2-\zeta)^{-1}\,
\right\|_{E_C}
\no \\
&\quad\qquad\leq\frac{\|B\|\|D\|_{E_C}}{({d}-\|B\|r)^2}
\|X_2-X_1\|,
\label{EEST}\quad X_1,X_2\in\cO_r.
\end{align}
Inequality \eqref{QQFC} implies that $F$ maps the ball $\cO_r$ into itself
whenever
\begin{equation}
\label{Est1}
\frac{\|D\|_{E_C}}{{d}-\|B\|r}\leq r.
\end{equation}
By \eqref{EEST} it follows that $F$ is a strict contraction on
$\cO_r$ if
\begin{equation}
\label{Est2}
\frac{\|B\|\|D\|_{E_C}}{({d}-\|B\|r)^2}<1\,.
\end{equation}
Solving \eqref{Est1} and \eqref{Est2} with respect to $r$
simultaneously, one obtains that if the radius of the ball
$\cO_r$ is within the bounds
\begin{equation}
\label{Br}
\frac{1}{\|B\|}\left(\frac{{d}}{2}-\sqrt{\frac{{d}^2}{4}-\|B\|\,\|D\|_{E_C}}
\right)
\leq  r < \frac{1}{\|B\|}
\left({d}-\sqrt{\|B\|\,\|D\|_{E_C}}
\right),
\end{equation}
then $F$ is a strictly contracting mapping of the ball $\cO_r$
into itself.  Applying Banach's Fixed Point Theorem one then infers
that equation (\Ref{EqQ}) has a unique solution in any ball
$\cO_r$ whenever $r$ satisfies \eqref{Br}.  Therefore, the fixed
point does not depend on the radii satisfying~(\Ref{Br}) and
hence it belongs to the smallest of these balls. This
observation proves the estimate
\begin{equation}
\label{QnEst}
\|X\|\leq\frac{1}{\|B\|}\left(
\frac{{d}}{2}-\sqrt{\frac{{d}^2}{4}-\|B\|\,\|D\|_{E_C}}\right).
\end{equation}
Finally, using
(\Ref{QQFC}),  for the fixed point  $X$ one
obtains the estimate
\begin{equation}
\label{QEC1}
\|X\|_{E_C}=\|F(X)\|_{E_C}\leq
\frac{\|D\|_{E_C}}{{d}-\|B\|\,\|X\|}.
\end{equation}
Then~(\Ref{QnEst}) yields
$$
\|X\|_{E_C}\leq\frac{\|D\|_{E_C}}
{\frac{{d}}{2}+\sqrt{\frac{{d}^2}{4}-\|B\|\,\|D\|_{E_C}}}=
\frac{1}{\|B\|}
\left({\frac{{d}}{2}-\sqrt{\frac{{d}^2}{4}-\|B\|\,\|D\|_{E_C}}}\right),
$$
proving \eqref{QEst1}.

Finally, using \eqref{BDest} and \eqref{QEst1} one easily verifies by
inspection that any of the assumptions \eqref{BDestSum} and
\eqref{BDestSum1} implies \eqref{Xstrict}.

The proof is complete.
\end{proof}
%%%%%%%%%%%%%%%%%%%%%%%%%%%%%%%%%%%%%%%%%%%%%%%%%%%%
\begin{remark}
This theorem extends results obtained in
\cite{AMM,MM99,MotSPbWorkshop,MotRem} where the operators $A$ and
$C$ are assumed to be self-adjoint. In case where the self-adjoint
operator $C$ is bounded, $D$ is a Hilbert-Schmidt operator, $B$ is
bounded, and $A$ is possibly unbounded densely defined closed
non-self-adjoint operator, the solvability of the equation
\eqref{EqQ1} under condition \eqref{tsemibd} has also been studied
in \cite{AdLT}.
\end{remark}

%%%%%%%%%%%%%%%%%%%%%%%%%%%%%%%%%%%%%%%%%%%%%%%%%%%%%%%%%%%
\vspace*{2mm} \noindent {\bf Acknowledgments.} This was supported
by the Deutsche For\-sch\-ungs\-gemeinschaft (DFG), the
Heisenberg-Landau Program, and the Russian Foundation for Basic
Research.  A. K. Motovilov gratefully acknowledges the kind
hospitality of the Institut f\"ur Angewandte Mathematik,
Universit\"at Bonn.

%%%%%%%%%%%%%%%%%%%%%%%%%%%%%%%%%%%%%%%%%%%%%%%%%%%%%%%

\end{document}